\documentclass{article} 

\usepackage{scalefnt}
\usepackage{amsmath,amssymb}
\usepackage{bbm}
\usepackage{wrapfig}
\usepackage{verbatim}
\usepackage{epsfig}
\usepackage{textcomp}
\usepackage{vmargin}
\usepackage{url}
\usepackage{xr}
\usepackage{cleveref}
\usepackage{algorithm}
\usepackage{algorithmic}
\usepackage{amsthm}

\usepackage{graphicx}          
 
\def\R{\mathbb{R}}

\def\E{\mathbb{E}}
\def\P{\mathbb{P}}
\def\A{\mathbb{A}}
\def\O{\mathbb{O}}
\def\X{\mathbb{X}}
\def\K{\mathbb{K}}
\def\T{\mathbb{T}}

\def\cd{\check{d}}

\newcommand{\thetabar}{{\bar{\theta}}}

\newcommand{\Psibar}{{\bar{\Psi}}}

\newcommand{\vbar}{{\bar{v}}}

\def\Mtn1{M_{t_{n+1}}}
\def\Mtn{M_{t_{n}}}
\def\Xtn1{X_{t_{n+1}}}
\def\Xtn{X_{t_{n}}}
\def\Ytn1{Y_{t_{n+1}}}
\def\Ytn{Y_{t_{n}}}

\def\phi0lip{[\Phi_0]}

\newcommand{\1}{\mathbbm{1}}
\newcommand{\cutline}{\qquad \qquad \qquad}

\usepackage[]{xcolor}
\definecolor{darkred}{rgb}{0.65,0.15,0.25}
\definecolor{darkgreen}{rgb}{0,0.4,0}
\definecolor{myblue}{rgb}{0.2,0.2,0.7}
\definecolor{aliceblue}{rgb}{0.94, 0.97, 1.0}
\definecolor{airforceblue}{rgb}{0.36, 0.54, 0.66}
\definecolor{antiquebrass}{rgb}{0.8, 0.58, 0.46}
\definecolor{cambridgeblue}{rgb}{0.64, 0.76, 0.68}
\definecolor{asparagus}{rgb}{0.53, 0.66, 0.42}
\definecolor{ceruleanblue}{rgb}{0.16, 0.32, 0.75}
\definecolor{chr1}{HTML}{A6CEE3}
\definecolor{chr2}{HTML}{1F78B4}
\definecolor{chr3}{HTML}{B2DF8A}
\definecolor{chr4}{HTML}{33A02C}
\definecolor{chr5}{HTML}{FB9A99}
\definecolor{chr6}{HTML}{E31A1C}
\definecolor{chr7}{HTML}{FDBF6F}
\definecolor{chr8}{HTML}{FF7F00}
\definecolor{chr9}{HTML}{CAB2D6}
\definecolor{chr10}{HTML}{6A3D9A}
\definecolor{chr11}{HTML}{FFFF99}
\definecolor{chr12}{HTML}{B15928}
\definecolor{chr13}{HTML}{8DD3C7}
\definecolor{chr14}{HTML}{FFFFB3}
\definecolor{chr15}{HTML}{BEBADA}
\definecolor{chr16}{HTML}{FB8072}
\definecolor{chr17}{HTML}{80B1D3}
\definecolor{chr18}{HTML}{FDB462}
\definecolor{chr19}{HTML}{B3DE69}
\definecolor{chr20}{HTML}{FCCDE5}
\definecolor{chr21}{HTML}{D9D9D9}
\definecolor{chr22}{HTML}{BC80BD}
\definecolor{chrX}{HTML}{CCEBC5}
\definecolor{chrY}{HTML}{FFED6F}

\usepackage{tikz}
\usetikzlibrary{decorations.markings}
\usetikzlibrary{decorations.pathreplacing}
\usetikzlibrary{spy}
\usetikzlibrary{positioning}
\usetikzlibrary{backgrounds}
\usetikzlibrary{hobby,decorations}
\usetikzlibrary{shapes,snakes}
\usetikzlibrary{shapes.geometric}
\usetikzlibrary{shadows}
\usetikzlibrary{intersections}
\usetikzlibrary{matrix}
\tikzset{   invisible/.style={opacity=0,text opacity=0},     visible on/.style={alt={#1{}{invisible}}},     alt/.code args={<#1>#2#3}{%
     \alt<#1>{\pgfkeysalso{#2}}{\pgfkeysalso{#3}} },}
\tikzset{cross/.style={cross out, draw=black, minimum size=2*(#1-\pgflinewidth), inner sep=0pt, outer sep=0pt},cross/.default={1pt}}

\def\newblock{\hskip .11em plus .33em minus .07em}

\newtheorem{theo}{Theorem}[section]
\newtheorem{lem}[theo]{Lemma}
\newtheorem{prop}[theo]{Proposition}

\newtheorem{hyp}[theo]{Assumption}

\newcommand{\bs}[1]{\boldsymbol{#1}}

\DeclareMathAlphabet\mathbfcal{OMS}{cmsy}{b}{n}
\DeclareMathOperator*{\argmin}{argmin}

\begin{document}

\title{Numerical method for solving impulse control problems in partially observed piecewise deterministic Markov processes\thanks{This work was partially supported by ANR Project HSSM-INCA (ANR-21-CE40-0005), by a European Union’s Horizon 2020 research and innovation program (Marie Sklodowska-Curie grant agreement No 890462)  and has been realized with the support of MESO@LR-Platform at the University of Montpellier.}}                                              

\author{Alice Cleynen \and Beno\^\i te de Saporta}       
\date{\small John Curtin School of Medical Research, The Australian National University, Canberra, ACT, Australia and IMAG, Univ Montpellier, CNRS, Montpellier, France}

\maketitle                       

\begin{abstract}
Designing efficient and rigorous numerical methods for sequential decision-making under uncertainty is a difficult problem that arises in many applications frameworks. In this paper we focus on the numerical solution of a subclass of impulse control problem for piecewise deterministic Markov process (PDMP) when the jump times are hidden. We first state the problem as a partially observed Markov decision process (POMDP) on a continuous state space and with controlled transition kernels corresponding to some specific skeleton chains of the PDMP. Then we proceed to build a numerically tractable approximation of the POMDP by tailor-made discretizations of the state spaces. The main difficulty in evaluating the discretization error comes from the possible random 
jumps of the PDMP between consecutive epochs of the POMDP and requires special care. Finally we 
discuss the practical construction of discretization grids and illustrate our method on simulations.
\end{abstract}

\newpage
\tableofcontents
\newpage
%
\section{Introduction}\label{sec:intro}
%
A large number of problems in science, including resource management, financial portfolio management, monitoring complex systems or medical treatment design to name just a few, can be characterized as sequential decision-making problems under uncertainty. In such problems, an agent interacts with a dynamic, stochastic, and incompletely known process, with the goal of finding an action-selection strategy that optimises some performance measure over several time steps. In optimal stopping problems, for instance when the agent has to decide when to replace some component of a production chain before full deterioration, the policy typically does not influence the underlying process until replacement when the process starts again with the same dynamics. However in many decision-making problems an important aspect is the effect of the agent's policy on the data collection; different policies naturally yielding different behaviors of the process at hand. 

In this paper we focus on the numerical solution of a subclass of impulse control problem for Piecewise Deterministic Markov Process (PDMP) when the process is partially observed, in the very challenging case when jump times are hidden. PDMPs are continuous time processes with hybrid state space: they can have both discrete and Euclidean variables. The dynamics of a PDMP is determined by three local characteristics: the flow, the jump intensity and the Markov jump kernel, see \cite[Section 24]{Davis84}. Between jumps, trajectories follow the deterministic flow.
The frequency of jumps is determined by the intensity, and the post-jump location is selected by the Markov kernel. 

General impulse control for PDMPs allows the controller to act on the process dynamics by choosing intervention dates and new states to restart from at each intervention. This family of problems was first studied by Costa and Davis in \cite{CD89}. It  has received a lot of attention since, see e.g.~\cite{Alm01,CR00,dSDZ16,Dempster95}, and was further extended for instance in \cite{DHP16,PPTZ19}. In these papers, the authors define a rigorous mathematical framework to state such control problem and establish some optimality equations such as dynamic programming equations for the value function. Numerical methods to compute an approximation of the value function and an $\epsilon$-optimal policy are also briefly or extensively described for instance in \cite{CD89,Aut12,Aut16}. They rely on discretization of the state space, either with direct cartesian grids or with dynamic grids obtained from simulations of the inter-jump-time-post-jump-location discrete time Markov chain embedded in the PDMP.

In all the papers cited above, the process is supposed to be perfectly observed at all times. 
{However in many real-life applications, such as dynamic reliability or biological systems, continuous measurements may not be available, and measures may be corrupted by noise. 
The field of dynamic reliability is particularly suitable for controlled PDMP models, see e.g. \cite{cocozza2021markov,dSDZ16,devooght1997dynamic} or \cite[p96]{lange2017cost}. Consider a complex system with interacting components subject to degradation or failure. The controller aims to design a maintenance policy to optimize some suitable target, such as maximizing production while minimizing downtime or maintenance costs. The discrete variables correspond to the states of the different components, and the Euclidean variables are those that influence the process, for instance the age of the components, or some physical characteristics such as pressure, temperature, \ldots. It is natural in this context to suppose that the partial or total failure of some components may be hidden. The physical variables are measured with captors, or monitoring systems and observations may thus be noisy, see e.g. \cite{baysse2013maintenance,BBGPPS14,soleimani2021integration}.}
{Biological systems can also be successfully modeled by PDMPs: one can cite examples in neurosciences \cite{RT15,RTW12}, population dynamics and genomics \cite{costa_piecewise_2016,DHKL15,herbach_inferring_2017}, or medical treatment \cite{Aut18,pasin_controlling_2018} to name just a few. Continuous-time observation of biological systems is often impossible, especially in-vivo, and again measurements may be indirect and corrupted by noise.}

Designing efficient and mathematically sound approximation methods to solve continuous time and continuous state space impulse control problems under partial observation is very challenging, especially for processes with jumps such as PDMPs. Relevant literature is scarce. One can mention \cite{BL17} or \cite{BdSD13} (in the special easier case of optimal stopping) where the position of the process is observed through noise, but the jump times are still perfectly observed, so that the properties of the inter-jump-time-post-jump-location chain can still be fully exploited. Another related recent work is \cite{C20} where the author studies an optimal control problem for pure jump processes (corresponding to PDMPs with constant flow) under partial observations. However, they consider continuous control instead of impulse control and the observations are not corrupted by noise.

A first step toward solving the impulse control problem for PDMPs under partial and noisy observations with hidden jump times was made by the authors in \cite{Aut18} in the special easier case of optimal stopping. They have shown that when trajectories of the process can be simulated, a double discretization allows to approximate the value function with general error bounds, and provides a candidate policy with excellent performance.
In this paper, we address a more general class of impulse control problems under partial observations where decisions do influence both the data collection and the dynamics of the process. {More specifically, we focus on a sub-class of impulse control problems with four specificities. First, the lapse between interventions can only take a finite number of values. Second, observations are collected at intervention times and corrupted by noise. For instance, in the dynamic reliability framework these two specificities correspond to to collecting data only at inspection dates, and the delay until the next inspection is scheduled in a given range (one month, one year, \ldots). Third, interventions only act on the discrete variables. For instance, in the dynamic reliability framework, we do not allow a maintenance that would reset the Euclidean variables to some nominal value. At an intervention, the settings of the system or the environment can be punctually modified leading to a different set of local characteristics without changing the value of the Euclidean variables. Finally, the state-space of the controlled PDMP has no active boundary, so that the process is punctuated only by random jumps and not by boundary jumps.}  

{This setting is more general than the setting of continuous-time MDPs (CTMDPs). Indeed, in CTMDPs, the state space is supposed to be finite or countable, whereas we have a hybrid state space here. In CTMDPs, the process evolves in continuous-time, and actions are taken at each jump time. The jump times are determined by exponential clocks from the controlled transition rates matrix and are not directly controller-chosen, see for instance \cite[Section 2.2]{guo2009continuous}. In our setting of controlled PDMPs, natural jumps may occur between interventions, an intervention does not correspond to a natural jump of the process, and intervention dates are directly controller-chosen.} 

{In our particular PDMP setting, the impulse control problem can be stated as a partially observed Markov decision process (POMDP) and this is the first main contribution of this paper. This POMDP is in discrete time, with epochs corresponding to intervention dates. However the state space is not discrete, and becomes infinite dimensional when turned into the corresponding fully observed Markov decision process (MDP) for the filter on the belief space, see e.g. \cite[Section 5.3]{bauerle11}. While such MDPs have theoretical exact optimal solutions, they are numerically intractable. In addition, the filter process is not simulatable, even for fixed policies, preventing the use of direct simulation-based discretization procedures. Our second main contribution is to propose and prove the convergence of an algorithm to approximate the value fonction and explicitly build a policy close to optimality. Our approach is based on tailor-made discretizations of the state spaces of the Euclidean variable of the PDMP and then of the corresponding discrete filter. It takes into account the major difficulty that the combinatorics associated to the possible decisions is too large to explore the whole belief space through simulations. Our third main contribution is to discuss the practical construction of the discretization grids, which is no easy problem.} Our results are illustrated by simulations on an example of medical treatment optimization for cancer patients.

The paper is organized as follows. 
In \cref{sec:statement} we state our generic optimization problem and turn it into a POMPD.
In \cref{sec:strategy} we give our main assumptions and  our resolution strategy, detailing construction of the approximate value function and candidate policy. 
{In \cref{sec:examples}, we provide a real-life example of medical treatment optimization fitting our framework and briefly present other general classes of applications fitting our framework.}
Experimental results and discussion of the practical construction of discretization grids are provided in \cref{sec:simulation}. Finally \cref{sec:conclusion} provides a short conclusion.
The main proofs are postponed to the appendix as well as details and specifications of the model used for the simulation study together with further discussion on grids construction.
%
\section{Problem statement}
\label{sec:statement}
%
We start with specifying the special class of impulse control problem for PDMPs we are focusing on, then we show how our control problem can be expressed as a POMDP.
%
\subsection{Impulse control for hidden PDMPs}
\label{sec:PDMP}
%
Let us first define the class of controlled PDMPs we consider. 
Let $\mathcal{M}=L\times M$ be a two-dimensional finite set. We will call regimes or modes its elements. We use a product set to distinguish between modes in $L$ that can be controller chosen and modes in $M$ that cannot.
For all regime $(\ell,m)$ in $\mathcal{M}$, let $E_m^\ell$ be a bounded Borel subset of $\mathbb{R}^d$ endowed with a norm $\|\cdot\|$. Set $E=\{(\ell,m,\mathtt{x}), \ell\in L, m\in M, \mathtt{x} \in E_m^\ell\}$, and $E^\ell=\{(m,\mathtt{x}), m\in M, \mathtt{x} \in E_m^\ell\}$ for all $\ell\in L$.
A PDMP on the state space $E$ is determined by three local characteristics:
\begin{itemize}
\item the flow $\bs\Phi(x,t)=(\ell,\bs\Phi^\ell(m,\mathtt{x},t))=(\ell,m,\Phi_m^\ell(\mathtt{x},t))$ for all $x=(\ell,m,\mathtt{x})$ in $E$ and $t\geq 0$, where $\Phi_m^\ell: \mathbb{R}^d \times \mathbb{R}_+ \rightarrow \mathbb{R}^d$ is continuous and satisfies a semi-group property $\Phi_m^\ell(\cdot,t+s)=\Phi_m^\ell(\Phi_m^\ell(\cdot,t),s)$, for all $t,s\in \mathbb{R}_+$. It describes the deterministic trajectory between jumps. 
\item the jump intensity $\bs\lambda(x)=\bs{\lambda}^\ell(m,\mathtt{x})=\lambda_m^\ell(\mathtt{x})$  for all $x=(\ell,m,\mathtt{x})$ in $E$, where $\lambda_m^\ell: {E}_m^\ell\rightarrow \mathbb{R}^+$ is a measurable function such that for any $x=(\ell,m,\mathtt{x})$ in $E$, there exists $\epsilon >0$ such that
\begin{align*}
 \int_0^{\epsilon} \lambda_m^\ell(\Phi_m^\ell(\mathtt{x},t))dt < +\infty. 
\end{align*}
For all $x=(\ell,m,\mathtt{x})$ in $E$ and $t\geq 0$, set 
\begin{align*}
\bs\Lambda(x,t)&=\bs\Lambda^\ell(m,\mathtt{x},t)=\Lambda_m^\ell(\mathtt{x},t)=\int_0^{t} \lambda_m^\ell(\Phi_m^\ell(\mathtt{x},s)) ds.
\end{align*}
\item the Markov kernel $\bs Q$ on $(\mathcal{B}({E}),{E})$ represents the transition measure of the process and allows to select the new location after each jump. It satisfies for all $x \in {E}$, $\bs Q(\{x\}|x)=0$. We also write $Q_{m}^{\ell}(\cdot|\mathtt{x})=\bs Q(\cdot|x)=\bs Q^\ell(\cdot|m,\mathtt{x})$ for all $x=(\ell,m,\mathtt{x})\in E$ and add the additional constraint that $Q$ cannot change the value of $\ell$, as $\ell$ is intended to be controller-chosen, i.e. $Q_{m}^{\ell}$ sends $E^\ell$ onto itself.
\end{itemize}

The formal probabilistic apparatus necessary to precisely define controlled trajectories and to formally state the impulse control problem is rather cumbersome, and will not be used in the sequel. Therefore, for the sake of simplicity, we only present an informal description of the optimal control problem. The interested reader is referred to \cite{CD89}, \cite[section 54]{Davis93} or \cite{DHP16} for a formal setting. Our optimization problem will be rigorously stated as a POMDP in \cref{sec:POMDP}.

We consider a finite horizon problem. Let $H>0$ be the optimization horizon, $\delta>0$ a fixed minimal lapse such that $N=H/\delta$ is an integer. Note that $\delta$ is not supposed to be small. A general impulse strategy $\mathcal{S}=(\ell_n,r_n)_{0\leq n\leq N-1}$ is a sequence of non-anticipative $E$-valued random variables on a measurable space $(\Omega,\mathcal{F})$ and of non-anticipative intervention lapses. In this work, we only consider a subclass of strategies where $\ell_n$ takes values in $L$ and $r_n$ is a multiple of $\delta$: $r_n\in\mathbb{T}$ where $\mathbb{T}$ is a subset of $\delta^{1:N}=\{\delta,2\delta,\ldots, N\delta\}$ that contains $\delta$. This means that on the one hand, the controller can only act on the process by changing its regime, i.e. by selecting the local characteristics to be applied until the next intervention, and on the other hand the lapse between consecutive interventions belongs to the finite set $\mathbb{T}$.
The trajectories of the PDMP controlled by strategy $\mathcal{S}$ are constructed recursively between intervention dates $\tau_n$ (defined recursively by $\tau_0=0$ and $\tau_{n+1}=\tau_{n}+r_n$) as described in \cref{algo:SimuPDMP}. 
\begin{algorithm}[tb]
\caption{Simulation of a trajectory of a controlled PDMP between interventions $n$ and $n+1$ from state $x_n=(\ell_n, m_n, \mathtt{x}_n)$}
\label{algo:SimuPDMP}
\begin{algorithmic}[1]
    \STATE $s\leftarrow 0$
    \STATE $\ell \leftarrow \ell_n, m\leftarrow m_n, \mathtt{x}\leftarrow \mathtt{x}_n, x\leftarrow (\ell,m,\mathtt{x})$
    \WHILE{$s<r_n$}
        \STATE $S\sim \lambda_m^{\ell}(\mathtt{x})$
        \STATE $X_{\tau_{n}+s+t} \leftarrow \bs\Phi(x,t)$ for $0\leq t\leq \min\{S,r_n-s\}$
        \STATE $\mathtt{x}\leftarrow \Phi_m^\ell(\mathtt{x},\min\{S,r_n-s\})$
        \IF{$S< r_n-s$}
            \STATE $s\leftarrow s+S$
            \STATE $(\ell,m',\mathtt{x}')\sim Q_m^\ell(\cdot | \mathtt{x})$
            \STATE $m\leftarrow m', \mathtt{x}\leftarrow \mathtt{x}', x \leftarrow (\ell,m,\mathtt{x})$
        \ELSE
        		\STATE $s\leftarrow r_n$
        \ENDIF
    \ENDWHILE
\end{algorithmic}
\end{algorithm}
%
In Line 5 of \cref{algo:SimuPDMP}, $S\sim \lambda_m^{\ell}(\mathtt{x})$  means that $S$ has the survival function
\begin{align*}
\mathbb{P}_x(S>t)=e^{-\int_0^t \lambda_m^{\ell}(\Phi_m^{\ell}(\mathtt{x},s))ds}=e^{-\bs\Lambda(x,t)}.
\end{align*}
{Note that boundary jumps are not allowed.}
The strategy $\mathcal{S}$ induces a family of probability measures $\mathbb{P}^{\mathcal{S}}_x$, $x\in E$, on a suitable probability space $(\Omega,\mathcal{F})$. Associated to strategy $\mathcal{S}$, we define the following expected total cost for a process starting at $x\in E$
\begin{align}\label{cout_strategie}
\mathcal{J}(\mathcal{S},x)=\E_x^{\mathcal{S}} \Bigl[ \int_0^{H} c_r(X_s)ds + \sum_{n=0}^{N-1}  c_i({X}_{\tau_n},\ell_n,r_n)+c_t(X_H)\Bigr],
\end{align}
where $\E_x^{\mathcal{S}}$ is the expectation with respect to $\mathbb{P}^{\mathcal{S}}_x$, $c_r$ is the running cost, $c_i$ the intervention cost and $c_t$ the terminal cost. 

The last ingredient needed to state the optimisation problem is to define admissible strategies. Again, the rigorous definition will be given in \cref{sec:POMDP} in the framework of POMDPs. Informally, decisions can only be taken in view of some discrete-time noisy observations of the process, instead of the exact value of the process at all times. More specifically, we assume that
\begin{itemize}
\item observations are only available at decision times $\tau_n$,
\item the controller-chosen regimes $\ell\in L$ are observed, the uncontrolled regimes $m\in M$ are hidden, 
\item the Euclidean variable is observed through noise: at time $\tau_n$, if the state of the process is $x=(\ell,m,\mathtt{x})$, the controller receives observation $y_n=F(x)+\epsilon_n=F(\mathtt{x})+\epsilon_n$ where $(\epsilon_n)$ are real-valued independent and identically distributed random variables with density $f$ independent from the controlled PDMP. We further assume that the random variables $y_n$ take values in a compact interval $I$ of the real line.
\end{itemize}


Denote by $\mathbb{S}$ the set of all admissible strategies. Our aim is to compute an approximation of the value function
\begin{align*}\label{valeur}
\mathcal{V}(x)=\inf_{\mathcal{S}\in\mathbb{S}}\mathcal{J}(\mathcal{S},x),
\end{align*}
and explicitly construct a strategy close to optimality. We also assume that the process starts from a fixed point $x_0=(\ell_0,m_0,\mathtt{x}_0)\in E$.
%
\subsection{Partially observed Markov decision process}
\label{sec:POMDP}
%
Finding a suitable rigorous way to state an impulse control problem for PDMPs with hidden jumps and under noisy observation is by no means straightforward, especially as regards defining admissible strategies, see e.g. \cite{Alm01} or \cite[sec. 1.1]{CD13}. This is our first main contribution in this paper.
As decision dates are discrete, we use the framework of POMDPs to rigorously state our control problem. In the sequel, we drop the regime $\ell\in L$ from the state $x$ and include it in the action instead to better fit the standard POMDP notation.
Denote $E_m=\cup_{\ell\in L}E_m^\ell$, 
and $E_M=\{(m,{\tt x}), m\in M, {\tt x} \in E_m\}$. 
Define on $E_M$ the following hybrid distance, for $x=(m,{\tt x})$ and $x'=(m',{\tt x}')\in E_M$,
\begin{align*}
\|x-x'\| = (\|{\tt x}-{\tt x}'\|) \1_{m=m'}+\infty \1_{m \neq m'}.
\end{align*}
%
Let $(\mathbb X,\mathbb A, \mathbb K, R,c,C)$ be the POMPD with the following characteristics.

$\bullet$ The state space is 
{$\mathbb X = \{\xi=(x,y,w)\in (E_M\times \O)\}\cup \{\Delta\}$}, with 
{$\O=\{\gamma \in I\times\delta^{1:N}\}$ the observation space}. It gathers the values $X_n$ of the hidden PDMP at dates $\tau_n$ as well as the observation process $Y_n$, and 
one additional observed variable 
$W_n$ 
that is the time elapsed since the beginning. State $\Delta$ is a cemetery state where the process is sent after the horizon $H$ 
is reached.

$\bullet$ The action space is $\mathbb A = (L\times\mathbb T)\cup \{\cd\}$, {where $\check{d}$ is an empty decision that is taken when the horizon $H$  
is reached. This purely technical decision sends the process to the cemetery state $\Delta$}.

$\bullet$ The constraints set $\mathbb K\subset \mathbb X\times \mathbb A$ is such that its sections $\mathbb K(\xi)=\{d\in \mathbb A; (\xi,d)\in\mathbb K\}$ satisfy $\mathbb K(\Delta)=\{\cd\}$ and 
$\mathbb K(x,y,w)=\K(w)$ for 
$(x,y,w)\in\X-\{\Delta\}$ as contraints concern the remaining control time. More specifically, one has 
\begin{itemize}
\item[--] $\K(w)\subset L\times(\mathbb{T} \cap \{\delta,2\delta,\dots,N\delta-w\})$ if $w\leq \delta(N-1)$, to force the last intervention to occur exactly at the horizon time $H$, 
\item[--] $\K(N\delta)=
\cd$, no intervention is possible after the horizon 
has been reached.
\end{itemize}
Note that as $\delta\in\mathbb{T}$,  the set $\mathbb{T} \cap \{\delta,2\delta,\dots,N\delta-w\}$ is never empty.

$\bullet$ The controlled transition kernels $R$ are defined as follows: for any bounded measurable function $g$ on $\X$, any $\xi\in\X$ and $d\in\K(\xi)$, one has
{\begin{align*}
Rg(\xi,\cd)&=g(\Delta),\\
Rg(\xi,d)&=\int_{I}\int_{E} g(x',y',w+r)f(y'-F(x'))P(dx'|x,d)dy',
\end{align*} }
if $\xi=(x,y,w)$, $d=(\ell,r)\neq\cd$, and $Rg(\xi,d)=0$ otherwise, 
where $P(\cdot| x,d)$ is the distribution of $X_r$ conditionally to $X_0=x$ under regime $\ell$, if $d=(\ell,r)$. 
 {Note that since conditionnaly on $X_n$, $Y_n$ is independent of any past $X_k$ and $Y_k$, kernel $R$ does not depend on $y$.}

$\bullet$ The terminal cost function 
$C:\mathbb X\to \mathbb R_+$ satisfies $C(\Delta)=0$ and $C(x,y,w)=c_t(x)$ for $(x,y,w)\in \X-\{\Delta\}$. 

$\bullet$ The cost-per-stage function $c:\mathbb K\times\mathbb X\to \mathbb R_+$ satisfies $c(\Delta, \cd,\cdot)=0$, $c(\xi,\cd,\cdot)=C(\xi)$ if $\xi\neq \Delta$. 
Ideally, to match \cref{cout_strategie}, for $d=(\ell,r)$ one should choose $c$ as 
\begin{align*}
c(x,d)=\E_x^{\ell} \Bigl[ \int_0^{r} c_r(X_s)ds\Bigr] + c_i(x,d).
\end{align*}
However, the integral part has no simple analytical expression and we choose some simpler proxys instead, see \cref{sec:cost}.

$\bullet$ The optimisation horizon is finite and equals $N$.\\

Classically, the sets of observable histories are defined recursively by $H_0=\O\cup\{\Delta\}$ and $H_n=H_{n-1}\times\mathbb{A}\times (\O\cup\{\Delta\})$.
A decision rule at time $n$ is a measurable mapping $g_n:H_n\rightarrow\mathbb{A}$ such that $g_n(h_n)\in\mathbb{K}(\gamma_n)$ for all histories $h_n=(\gamma_0,d_0,\gamma_1,d_1,\ldots,\gamma_n)$. A sequence $\pi=(g_n)_{0:N-1}=(g_0,\ldots,g_{n},\ldots, g_{N-1})$ where $g_k$ is a decision rule at time $k$  is called an admissible policy. Let {$\Pi_N$} denote the set of all admissible policies. 
The controlled trajectory of the POMDP following policy $\pi=(g_n)_{0:N-1}\in\Pi_N$ is defined recursively by $\Xi_0\in \X$ and for $0\leq n\leq N-1$,
\begin{itemize}
\item $A_{n}=g_n(\Xi_n)$,
\item $\Xi_{n+1}\sim R(\cdot | \Xi_{n},A_{n})$.
\end{itemize}
Note that the cemetery state $\Delta$ ensures that all trajectories have the same length $N$, even if they do not have the same number of actual decisions ($d\neq \cd$).
Then one can define the total expected cost of policy $\pi\in\Pi_N$ starting at $\xi_0\in \X$ as
\begin{equation*}
 J_{\pi}(\xi_0)=\mathbb{E}^{\pi}_{\xi_0}\left[\sum_{n=0}^{N-1} c(\Xi_n,A_n,\Xi_{n+1})+C(\Xi_N)\right],
\end{equation*}
 and our control problem corresponds to the optimisation problem
\begin{equation*}
V(\xi_0)=\inf_{\pi\in \Pi_N}J_\pi(\xi_0).
\end{equation*}
%
Now the problem is rigorously stated, our next aim is now to compute an approximation of the value function $V$ 
and explicitly construct a strategy close to optimality. 
%
\section{Resolution strategy and assumptions}
\label{sec:strategy}
%
Our second main contribution is to propose a numerical approach to (approximately) solve our POMDP.
The first difficulty to solve our POMDP comes from the fact that it is partially observed. Thus our first step is to convert it into an equivalent fully observed MDP on a suitable belief space $\X'$ by introducing a filter process. This is done in \cref{sec:belief}.
While dynamic programming equations hold true for the fully observed MDP and in theory provide the exact optimal strategy to the optimisation problem, there are two main difficulties to its practical resolution. On the one hand, the state space $\X'$ is continuous and infinite-dimensional. And on the other hand the filter process is not simulatable, even for fixed policies, preventing the use of direct simulation-based discretization procedures.

Our approach to (approximately) solve the belief MDP is to construct a numerically trac\-ta\-ble approximation of its value function based on discretizations of the dynamic programming equations. 
To obtain a finite-dimensional process and a simulatable approximation of the filter process, we first discretize the Euclidean part of the state space $E$ of the PDMP $(X_t)$. 
Then we discretize the state space of the approximated filtered process in order to obtain a numerically tractable approximation. Finally, we obtain an explicit candidate policy by solving the resulting discretized approximation of the dynamic programming equations. This procedure is detailed in \cref{sec:proof} {and illustrated in \cref{fig:proof}.}

\begin{figure}[tp]
\begin{center}
     \begin{tikzpicture}[,scale=1, every node/.style={scale=0.8}]
{      \node(cont) at (-1, 3) {${X}_t=(m_t,\mathtt{x}_t)$};
      \node(in) at (-1,1.5) {${X}_{n}=(m_{t_n},\mathtt{x}_{t_n})$};
      \node(in-EP) at (-1, 1) {Space $E_M$, Kernel $P$};
      \draw[->] (cont)--(in) node[left=0.15cm,pos=0.5] {\small $r\in\T$} ;
      \node(va) at (-1,-1.5) {$({X}_{n},{Y}_n,W_n)$};
      \node(va-EP) at (-1, -2) {Space $\mathbb{X}= E_M\times\mathbb{O}\cup\{\Delta\}$, Kernel ${R}$};
      \draw[->] (in-EP) --(va) node[left=0.15cm,pos=0.5,color=blue!50] {observations} ;
      \node(zn) at (0.1, 0.2) {\small{${Y}_n=f(\mathtt{x}_{t_n})+\epsilon_n$}};
      \node(zn) at (-0.2, -0.6) {{\small{{$W_n=w_{t_n}$}}}};
      \node(fi) at (-1,-3.5) {$(\Theta_n,Y_n,W_n)$};
      \node(fi-EP) at (-1, -4) {Space $\mathbb{X}'\subset \mathcal{P}(E_M)\times\mathbb{O}\cup\{\Delta\}$, Kernel ${R}'$};
      \node(dyn) at (-1,-5.5) {$V'_n(\Theta_n,Y_n,W_n)$};
      \draw[->] (va-EP)--(fi) node[left=0.15cm,pos=0.5,color=blue!50] {\small filtering} ;
      \draw[->] (fi-EP)--(dyn) node[left=0.15cm,pos=0.35,color=blue!50] {\small dynamic} node[left=0.15cm,pos=0.65,color=blue!50] {\small programming}; 
 }     
{       
      \node(q1) at (4,1.5) {$\bar{X}_{n}=({m}_{t_n},\mathtt{x}_{t_n})$};
      \node(q1-EP) at (4,1) {Finite space $\Omega$, Kernel $\bar{{P}}$};%
      \node(vaq1) at (4,-1.5) {$(\bar{X}_{n},\bar Y_n,{W}_n)$};
      \node(vaq1-EP) at (4, -2) {$\bar\X= \Omega\times\mathbb{O}\cup\{\Delta\}$, Kernel $\bar{R}$};
      \node(fiq1) at (4,-3.5) {$(\bar{\Theta}_n,\bar Y_n,{W}_n)$};
       \node(fiq1-EP) at (4, -4) {$\bar\X'\subset \mathcal{P}(\Omega)\times\mathbb{O}\cup\{\Delta\}$, Kernel $\bar{R}'$};
      \node(dynq1) at (4,-5.5) {$\bar V'_n(\bar{\Theta}_n,\bar Y_n,{W}_n)$}; 
      \draw[->] (in-EP)--(q1-EP) node[below=0.1cm,pos=0.5,color=blue!50] {\small {discretization 1}}; 
      \draw[->, thick,dotted] (va-EP)--(vaq1-EP) ;
      \draw[->, thick,dotted] (fi-EP)--(fiq1-EP);
	\draw[->, thick,dotted] (dyn)--(dynq1);
	\draw[->, thick,dotted]  (q1-EP)--(vaq1);
	\draw[->, thick,dotted]  (vaq1-EP)--(fiq1);
	\draw[->, thick,dotted] (fiq1-EP)--(dynq1);
}	
{  	
      \node(fiq2) at (8,-3.5) {$(\hat{\Theta}_n,{\bar Y_n},W_n)$};
       \node(fiq2-EP) at (8, -4) {$\Gamma$, Kernel $\hat{R}'$};
      \node(dynq2) at (8,-5.5) {$\hat V'_n(\hat{\Theta}_n,{\bar Y_n},W_n)$};
      \draw[->] (fiq1-EP)--(fiq2-EP) node[below=0.15cm,pos=0.5,color=blue!50] {\small {discretization 2}};
      \draw[->, thick,dotted] (dynq1)--(dynq2) ;
      \draw[->, thick,dotted] (fiq2-EP)--(dynq2);
      \draw[->, thick,dotted] (dynq1)--(dynq2);
 }     
    \end{tikzpicture}
\end{center}
\caption{{{\bf Graphical sketch of the approximation procedure}. The first column corresponds to the original PDMP / POMDP / fully observed MDP. The second column corresponds to the first discretization. The third column corresponds to the second discretization.}}\label{fig:proof}
\end{figure}

In order to keep track of the discretization errors throughout the different steps described above, we need our transition kernels to be regular enough. Hence we start this section with stating the assumptions required on the parameters of the PDMP. 
%
%
\subsection{Regularity assumptions}
\label{sec:regularity}
%
{The first set of assumptions will ensure the regularity of kernel $P$. First, we require the local characteristics of the PDMP to be regular enough. This is a classical assumption, that usually ensures the regularity of the kernel of the inter-jump-time-post-jump-location Markov chain embedded in the PDMP, see e.g. \cite[Section 1.4]{dSDZ16}. However, here kernel $P$ is more complex as any number of jumps may occur between two interventions. Therefore we need an additional restrictive assumption (assumption \ref{H-jumps}) to limit the number of jumps between consecutive interventions. Denote $\overline\delta=\max\{r\in \T\}$, the maximum time between consecutive interventions. 
%
\begin{hyp}\label{H-Phi} 
The flow $\bs \Phi$ is Lipschitz continuous: there exists a positive constant $[\Phi]$ 
such that for all $\ell \in L$, $x,x'\in E$, $0\leq t, t'\leq \overline\delta$,  
one has
\begin{align*}
|\bs \Phi^\ell(x,t)-\bs\Phi^\ell(x',t')|&\leq {[\Phi](\|x-x'\|+|t-t'|)}.
\end{align*}
\end{hyp}
%
\begin{hyp}\label{H-lambda-lip}
The jump intensity $\bs \lambda$ is Lipschitz continuous along the flow: there exists a positive constant $[\lambda]$  such that for all $\ell\in L$, $x$ and $x'\in E_M$, $0\leq t\leq \overline\delta$, one has
\begin{align*}
|\bs \lambda^\ell(\Phi^\ell(x,t))-\bs \lambda^\ell(\Phi^\ell(x',t))|\leq [\lambda]\|x-x'\|.
\end{align*}
\end{hyp}
%
\begin{hyp}\label{H-lambda-b}
The jump intensity $\bs \lambda$ is bounded: there exists a positive constant  $\|\lambda\|$ such that for all  $x\in E_M$, and $\ell\in L$, one has
\begin{align*}
 |\bs \lambda^\ell(x)|\leq \|\lambda\|.
\end{align*}
\end{hyp}
%
In the sequel, for all bounded measurable function $h$ we note 
\begin{align*}
\bs Q^\ell h(x)=\int_{E_M}h(x')\bs Q^\ell(dx',x).
\end{align*}
%
\begin{hyp}\label{H-Q}
The Markov kernel $\bs Q$ is Lipschitz continuous along the flow in the following sense: there exists a positive constant $[Q]$ such that for all functions $h\in BL_1(E)$, for all $\ell\in L$,
for all $x,x'\in E_M$, one has
\begin{equation*}
\left|\bs Q^\ell h(x)-\bs Q^\ell h(x')\right|
\leq [Q]\|x-x'\|.
\end{equation*}  
\end{hyp}
%
\begin{hyp}\label{H-jumps}
For any control $\ell\in L$, for any $x\in E_M$, the maximum number of jumps of the PDMP starting at time $0$ from $x$ and following control $\ell$ up to time $\overline\delta$ is bounded by $N_j<+\infty$.
\end{hyp}
%
This last assumption is further discussed in \cref{app:reg-P}. } 
We also require additional  regularity assumptions on  the observation process and the cost functions.
\begin{hyp}\label{H1} 
There exist non negative real constants $L_Y$, $\underline{f}$ and $\overline{f}$ 
such that for all {$(x,x') \in E_M$} and $y\in I$ one has
\begin{align*}
  \big|f(y-F({x}))-f(y-F({x'}))\big| &\leq L_Y\|{x-x'}\|,\\
0 < \underline{f} \leq  f(y-F({x}))  &\leq  \overline{f}<+\infty.
\end{align*}
We also set $L_f=L_Y|I|$ and $B_f=\overline{f}|I|$, where $|I|=\int_I dy$ is the length of interval $I$.
\end{hyp}
%
\begin{hyp}\label{Hc} 
There exist non negative real constants $L_c$, $L_C$, $B_c$ and $B_C$ such that for all $x,x',x''$ in $E_M$ and $d\in L\times \mathbb{T}$, one has
\begin{align*}
|C(x)|&\leq B_C,\\
|c(x,d,x')|&\leq B_c,\\
|C(x)-C(x')|&\leq L_C\|x-x'\|,\\
|c(x,d,x')-c(x'',d,x')|&\leq L_c\|x-x''\|,\\
|c(x,d,x')-c(x,d,x'')|&\leq L_c\|x'-x''\|.
\end{align*}
\end{hyp}
\subsection{Equivalent MDP on the belief state}
\label{sec:belief}
%
To convert a POMDP into an equivalent fully observed MDP is classical therefore details are omitted. The interested reader may consult e.g.  \cite[Section 5.3]{bauerle11} or \cite{BdSD13,Aut18} for similar derivations.
For $n\leq N$, set 
$\mathcal{F}_n^O=\sigma(Y_k,W_k,0\leq k\leq n)$  
the $\sigma$-field generated by the observations up to~$n$. Let 
\begin{align*}
\Theta_n(A)&=\mathbb{P}(X_{n}\in A | \mathcal{F}_n^O) = \mathbb{P}(\Xi_{n}\in A\times\{(Y_n,W_n)\} | \mathcal{F}_n^O),
\end{align*}
denote the filter or belief process for the unobserved part of the process, for any Borel subset $A$ of $E_M$. The standard prediction-correction approach yields a recursive construction for the filter.
%
\begin{prop}\label{def:Psi}
{ For any $n\ge 0$, 
conditionally on $(Y_{n+1},W_{n+1})=(y',w')$, $d=(\ell,r)\in L\times\mathbb T$ and $\Theta_{n}=\theta$, one has 
$\Theta_{n+1}=\Psi(\theta,y',w',d)$ with
\begin{align*}
  \Psi(\theta,y',w',d)(A) 
  &= \frac{\int_{E_M} \int_{E_M}f(y'-F(x'))\1_A(x')P(dx'|x,d)\theta(dx)}{\int_{E_M} \int_{E_M}  f(y'-F(x'))P(dx'|x,d)\theta(dx)}, 
\end{align*}
for any Borel subset $A$ of $E_M$.
}
\end{prop}\label{correction}
%
Let  $\mathcal{P}(A)$ denote the set of probability measures on set $A$. 
The equivalent fully observed MDP $(\mathbb X',\mathbb A, \mathbb K', R',c',C')$ is defined as follows.

$\bullet$ {The state space is a subset  $\X'$ of $(\mathcal{P}(E_M)\times \O)\cup\{\Delta\}$ satisfying the following constraints: all $\xi=(\theta,y,w)\in\X'$ satisfy 
$\theta(E_{m_0})\geq \big({\underline f}{\overline f}^{-1}\big)^{\frac{w}{\delta}\vee 1}e^{-w\|\lambda\|}$, where $m_0$ is the mode of the fixed starting point. }
It is necessary to restrict the state space to ensure the regularity of our operators, see \cref{app:reg-R'}.
This constraint is technical, and is guaranteed as soon as the process starts in mode~$m_0$, see Appendix \ref{app:R'X'}.

$\bullet$ The action space is still $\mathbb A = (L\times\mathbb T)\cup \{\cd\}$.

$\bullet$ The constraints set $\mathbb K'\subset \mathbb X'\times \mathbb A$ is such that its sections $\mathbb K'(\xi)=\{d\in \mathbb A; (\xi,d)\in\mathbb K'\}$ satisfy $\mathbb K'(\Delta)=\cd$ and $\mathbb K'(\theta,y,w)=\K'(w)=\K(w)$ for $(\theta,y,w)\in \mathcal P(E_M)\times\O$.

$\bullet$  The controlled transition kernels $R'$ are defined as follows: for any bounded measurable function $g$ on $\X'$, any $\xi\in\X'$ and $d\in\K'(\xi)$, one has $R'g(\xi,d)=g(\Delta)$ if $d=\cd$, and
{
\begin{align*}
R'g(\xi,d)&=\int_{E_M}\int_{I}\int_{E_M} g\big(\Psi(\theta,y',w+r,d),y',w+r\big) 
f(y'\!-F(x'))P(dx'|x,d)dy'\theta(dx) 
\end{align*}
}
if $\xi=(\theta,y,w)$, $d=(\ell,r)\neq \cd$.
See \cref{app:R'X'} for the proof that $R'$ maps $\X'$ onto itself.  

$\bullet$  The non-negative cost-per-stage function $c':\mathbb K'\to \mathbb R_+$ and the terminal cost function $C':\mathbb X'\to \mathbb R_+$ are defined by  $C'(\Delta)=c'(\Delta, \cd,\cdot)=0$ and for $\xi=(\theta,\gamma)\in \mathcal P(E_M)\times \O$ and $d\in \mathbb K'(\xi)$,
\begin{eqnarray*}
  c'(\xi,d)=c'(\theta,\gamma,d)&=&\int_{E_M^2} c(x,d,x') P(dx'|x,d)\theta(dx), \\
  C'(\xi)=C'(\theta,\gamma)&=&   \int_{E_M} C(x) \theta(dx).
\end{eqnarray*}

$\bullet$  The optimisation horizon is still  $N$.

\noindent Denote $(\Xi'_n)$ a trajectory of the fully observed MDP.
The cost of 
strategy $\pi\in\Pi_N$ 
is
\begin{equation*}
  J'(\pi,\xi'_0) = {\mathbb{E}}^{\pi}_{\xi'_0}\left[\sum_{n=0}^{N-1} c'(\Xi'_{n},A_n)+C'(\Xi'_{N})\right],
\end{equation*}
and the value function of the fully observed problem is,
\begin{equation*}
V'(\xi'_0)=\inf_{\pi\in \Pi_N}J'(\pi,\xi'_0).
\end{equation*}
%
If $\xi'_0=(\delta_{x_0},y_0,w_0)$, one has $V'(\xi'_0)=V(x_0,y_0,w_0)$,
so that solving the partially observed MDP is equivalent to solving the fully observed one.
In addition, the value function $V'$ satisfies the well known dynamic programming equations, see e.g. \cite[Section 5.3]{BL17}.
%
\begin{theo}\label{th:DP'}
For $\xi\in\X'$, set $v_N'(\xi)= C'(\xi)$  and for $0\leq n\leq N-1$, define by backwards induction
  \begin{eqnarray*}
    v_n'(\xi)&=&\min_{d\in\mathbb{K'}(\xi)} \left\{ c'(\xi,d) + R' v'_{n+1}(\xi,d) \right\}.
  \end{eqnarray*}
Let $\xi_0'=(\delta_{(m_0,\mathtt{x}_0)},y_0,0) \in \X'$. Then we have $v'_0(\xi_0') = V'(\xi_0')=V((m_0,\mathtt{x}_0),y_0,0).$ 
\end{theo}
%
\subsection{Approximations of the value function}
\label{sec:proof}
%
Our approximation of the value function is based on discretizations of the underlying state spaces. 
%
%
\subsubsection{First discretization.}
\label{sec:first}
%
Let $\Omega=\{\omega^1,\ldots,\omega^{n_{\Omega}}\}$ be a finite grid on $E_M$ containing at least one point in each mode $m\in M$. Let $p_\Omega$ denote the nearest-neighbor projection from $E_M$ onto $\Omega$ for the distance defined in \cref{sec:PDMP}. 
In particular, $p_\Omega$ preserves the mode. Set $\Omega_m=E_m\cap\Omega$ for all $n\in M$. Let $(C_i)_{1\leq i\leq K}$ be a Voronoi tessellation of $E_M$ associated to $\Omega$. Namely, $(C_i)_{1\leq i\leq n_{\Omega}}$ is a partition of $E_M$ such that for all $1\leq i\leq n_{\Omega}$, one has
\begin{equation*} 
    C_i\subset \{x\in E_M ; \|x-\omega^i\|\leq \|x-\omega^j\| \ \forall j\in 1:n_\Omega\}.
\end{equation*}
%
We will denote $\mathcal{D}_i$ the diameter of cell $C_i$: $\mathcal{D}_i=\sup \{\|x-x'\| ; x,x' \in C_i \}$.
%
{As stated in \cref{sec:PDMP}, the state spaces are supposed to be bounded. Denote $\|E_M\|=\sup\{\|x\|, x\in E_M\}$.
Hence all $\mathcal{D}_i$ are finite and bounded by $\|E_M\|$. }
%
%
%

We define the controlled kernels $\bar P$ from $E_M\times (L\times \mathbb{T})$ onto $\Omega$ as
\begin{equation*}
    \bar P(\omega^j | x,d) =P(C_j|x,d),
\end{equation*}
for all $x\in E$, $d\in L\times \mathbb{T}$ and $1\leq j\leq n_{\Omega}$. In particular, the restriction of $\bar P$ to $\Omega$ is a controlled Markov kernel on $\Omega$. We now replace kernel $P$ by kernel $\bar P$ in the dynamic programming equations of Theorem \ref{th:DP'} in order to obtain our first approximation.

{Let $\bar\X'\subset\{\xi=(\theta,\gamma)\in \X'; \theta(\Omega)=1, \}\cup\{\Delta\}$ satisfying the following constraints: all $\bar\xi=(\bar\theta,y,w)\in\bar\X'$ satisfy 
$\bar\theta(\Omega_{m_0})\geq \Big({\underline f}{\overline f}^{-1}\Big)^{\frac{w}{\delta}\vee 1}e^{-w\|\lambda\|}$, where $m_0$ is the mode of the fixed starting point. }
Let $\bar\Psi$ be the approximate filter operator defined by replacing the integrals w.r.t. $P$ in the filter operator $\Psi$ defined in Proposition \ref{def:Psi} by integrals w.r.t. $\bar P$. 
The approximated controlled transition kernels $\bar R'$ are defined as follows: for any bounded measurable function $g$ on $\bar \X'$, any $\xi\in\bar\X'$ and $d\in\K'(\xi)$, one has  $\bar R'g(\xi,d)=g(\Delta)$ if $d=\cd$, and
{
\begin{align*}
\bar R'g(\xi,d)&=\int_{I}\sum_{\omega^i\in\Omega}\sum_{\omega^j\in\Omega} g\big(\Psibar(\thetabar,y',w+r,d),y',w+r\big)f(y'-F(\omega^j))\bar P(\omega^j|\omega^i,d)\thetabar(\omega^i)dy',
\end{align*}
}%
if $\xi=(\thetabar,y,w)$, $d=(\ell,r)$.
Note that this kernel does not depend on $y$, and see \cref{app:R'X'bar} for the proof that $\bar R$ sends $\bar \X'$ onto itself.
The approximated cost-per-stage function is defined by $\bar c'(\Delta,\cd)=0$ and 
\begin{align*}
  \bar c'(\thetabar,\gamma,d)&=\sum_{\omega^i\in\Omega} \sum_{\omega^j\in\Omega}c(\omega^i,d,\omega^j) \bar P(\omega^j|\omega^i,d)\thetabar(\omega^i), 
\end{align*}
for all $\xi=(\thetabar,\gamma)\in\bar\X'$ and $d\in\K'(\xi)$.
Finally, for all $\xi\in\bar \X'$, set $\bar v_N'(\xi)=  \ C'(\xi)$,
 and for $0\leq n\leq N-1$, define by backwards induction
 \begin{eqnarray*}
   \bar v_n'(\xi)&=&\min_{d\in\mathbb{K}(\xi)} \left\{ \bar c'(\xi,d) + \bar R' \bar v'_{n+1}(\xi,d) \right\}.
 \end{eqnarray*}
%
If the grid $\Omega$ is precise enough, $\bar P$ should be a good approximation of kernel $P$ and thus one can expect that functions $\bar v'_n$ are good approximations of our value functions $v'_n$. Indeed, we have the following result.
%
\begin{theo}\label{theorem1}
Under Assumptions {\ref{H-Phi} to \ref{Hc}}, for all $\xi'\in\bar\X'$, one has
  \begin{align*}
   |v'_N(\xi') - \vbar_N'(\xi')|   & = 0 \\
   |v'_n(\xi')- \vbar_n'(\xi')| &\leq  C_{v_n'}\sup_{j \in \{1,\dots,n_\Omega\}} \mathcal{D}_j, \quad 0\leq n < N,
\end{align*}
where $C_{v'_n}$ depends only on $n$, $N$, $\delta$ and the regularity constants of the parameters.
\end{theo}
%
 Its proof is given in \cref{app:proof-th1}
 and is based on {the form of the kernels, their regularity (especially that of kernel $P$) and the dynamic programming equation}.
 
The main gain with this first approximation is that the filter operator $\bar \Psi$ now involves only finite weighted sums and therefore the corresponding approximate filter process is now simulatable. However, functions $\bar v'_n$ still cannot be computed because of the continuous integration in $y'$ and because the state space $\bar\X'$ is still continuous. We now proceed to a second discretization in order to obtain numerically tractable approximations of our value functions.
%
\subsubsection{Second discretization.}
\label{sec:second}
%
Let $\Gamma=\{\rho^1,\ldots,\rho^{n_{\Gamma}}\}$ be a finite grid on $\bar \X'$ containing at least $\Delta$, one element $\xi=(\theta,\gamma)$ such that $\theta(E_M)=1$. 
Let $p_\Gamma$ denote the nearest-neighbor projection from $\bar \X'$ onto $\Gamma$. Here, the distance on $\bar \X'$ is such that
$d(\xi, \xi')=d_\Omega(\thetabar,\thetabar')+|y-y'|$ if  $\xi=(\thetabar,y,w),\ \xi'=(\thetabar',y',w')$, with $w=w'$, $d(\xi,\xi')= +\infty$ otherwise,
where $d_\Omega$ is defined as
\begin{align*} 
d_\Omega(\thetabar,\thetabar')&=\left[\sum_{k=1}^{n_{\Omega}}|\thetabar(\omega^k)-\thetabar'(\omega^k)|^2\right]^{1/2},
\end{align*}
meaning that $\mathcal P(\Omega)$ is identified to the ${n_{\Omega}}$-dimensional simplex endowed with the $L^2$-norm.
Let $(\bar C_i)_{1\leq i\leq {n_{\Gamma}}}$ be a Voronoi tessellation of $\bar \X'$ associated to the grid $\Gamma$. We will denote $\mathcal{\bar D}_i$ the diameter of cell $\bar C_i$.
Note that the space $\bar \X'$ is compact, so that all $\mathcal{\bar D}_i$ are finite and bounded by $1+|I|$.
The controlled kernels $\hat R'$ are defined on $\bar \X'$ as follows: for any $\rho^j\in\Gamma$, any $\xi\in\bar\X'$ and $d\in\K'(\xi)$, one has
\begin{equation*}
    \hat R'(\rho^j | \xi,d) = \bar R'(\bar C_j|\xi,d).
\end{equation*}
In particular, the restriction of $\hat R'$ to $\Gamma$ is a family of controlled Markov kernels on $\Gamma$. For $1\leq j\leq {n_{\Gamma}}$, set $ \hat v_N'(\rho^j)=  C'(\rho^j)$,
 and for $0\leq n\leq N-1$ and $1\leq j\leq {n_{\Gamma}}$, define by backwards induction
 \begin{eqnarray*}
   \hat v_n'(\rho^j)&=&\min_{d\in\mathbb{K}'(\rho^j)} \left\{ \bar c'(\rho^j,d) + \hat R' \hat v'_{n+1}(\rho^j,d) \right\}.
 \end{eqnarray*}
%
Again, if the grid $\Gamma$ is precise enough, $\hat R'$ should be a good approximation of kernel $\bar R'$ and thus one can expect that functions $\hat v'_n$ are good approximations of functions $\bar v'_n$ and thus of our value functions $v'_n$. More precisely, we have the following result.
%
\begin{theo}\label{theorem2}
Under Assumptions \ref{H1} and \ref{Hc}, for all points $\rho^j\in \Gamma$, one has
\begin{align*}
   |\hat v'_N(\rho^j) - \vbar_N'(\rho^j)|   & =0, \\
   |\hat v'_n(\rho^j)- \vbar_n'(\rho^j)| &\leq C_{\bar v_n'}\sup_{j \in \{1,\dots,{n_{\Gamma}}\}} \mathcal{\bar D}_j, 
\end{align*}
where $C_{\bar v'_n}$ depends only on $n$, $N$, $\delta$ and the regularity constants of the parameters.
\end{theo}
%
 Its proof is given in \cref{app:proof-th2}.
The main gain with this second approximation is that integration against kernel $\hat R'$ boils down to computing finite weighted sum, and functions $\hat v'_n$ are defined on a finite state space which makes the dynamic programming equations fully tractable numerically. 
%
\subsection{Candidate strategy}
\label{sec:policy}
%
We can now construct a computable strategy using the fully discretized value function. The idea is to first build an approximate filter using the operator $\bar \Psi$ and project the resulting filter together with the current observation onto grid $\Gamma$. Then one selects the next decision using the dynamic programming equation on the grid $\Gamma$. More precisely, 
suppose that the process starts from point $\xi_0=(m_0,\mathtt{x}_0,y_0,0)$, such that the initial observation is $\gamma_0=(y_0,0)$.
One can recursively compute an approximate filter $(\bar\theta_n)$ and the corresponding decisions $(d_n)$ as follows.
First, set
\begin{align*}
\bar\theta_0=\delta_{(m_0,\mathtt{x}_0)},\quad
d_0(\bar \theta_0,\gamma_0)=\argmin_{d\in\K'(\bar\theta_0,\gamma_0)} \big\{\bar c'(p_{\Gamma}(\bar \theta_0,\gamma_0),d)+\hat R'\hat v'_{1}(p_{\Gamma}(\bar \theta_0,\gamma_0),d) \big\}.
\end{align*}
Suppose one has constructed the sequence $(\bar\theta_n, d_n)$ up to stage $k-1$. Then after receiving the $k$-th observation $\gamma_k$, the next approximate filter and decisions are
\begin{align*}
\bar\theta_k=\bar\Psi(\bar\theta_{k-1},\gamma_k,d_{k-1}),\ 
d_k(\bar \theta_k,\gamma_k)=\!\!\!\argmin_{d\in\K'(\bar\theta_0,\gamma_0)}\!\!\! \big\{\bar c'(p_{\Gamma}(\bar \theta_k,\gamma_k),d)+\hat R'\hat v'_{1}(p_{\Gamma}(\bar \theta_k,\gamma_k),d) \big\}.
\end{align*}
with the convention that the last (not required) decision is $d_N=\cd$.
Note that it should be better to use operator $\bar\Psi$ on $\bar\theta_{k-1}$ then project the result onto grid $\Gamma$ than using operator $\hat \Psi$ (obtained by replacing integration wrt $\theta$ by integration wrt $p_\Gamma(\theta)$ in the definition of $\bar\Psi$) as it should generate a smaller error. Although it is reasonable to think that this strategy should be close to optimality, it is  an open problem to actually prove it as the sequence $(\bar\theta_n,\gamma_n)$ is not generated by the kernel $\bar R'$. Indeed, we use here the observations generated by the original sequence $X_n$ with kernel $P$ and not that generated by kernel $\bar P$, as in practice only the original observations are available. Its performance is assessed in \cref{sec:simulation} on a specific example. 
%
\section{Examples of impulse control problems}
\label{sec:examples}
%
In this section, we give define the specific medical example on which we ran our approximation procedure, and give some details of other control problems fitting into our framework.
%
\subsection{Optimisation of patient follow-up.}
\label{sec:medexample}

We consider an example of patient follow-up that we will use in the Simulation Study. The
mode $m$ corresponds to the overall state of the patient ($m=0$: sound, $m=1$: disease 1, $m=2$: disease 2, $m=3$: death of the patient). The Euclidean variable is two-dimensional $\mathtt{x}=(\zeta,u)$. Variable $\zeta\in[K_i,K_s]$ corresponds to some marker of the disease that can be measured, and $u\geq 0$ the time since the last jump in order to encompass semi-Markov dynamics. The control $\ell$ correspond to the medical treatment ($\ell=\emptyset$: no treatment, $\ell=a$: efficient for disease $1$ and {slows the progression of disease $2$}, $\ell=b$: efficient for disease $2$ and {slows the progression of disease $1$). Decision dates correspond to visits to the medical center when the marker is measured and a new treatment is selected and applied until the next visit. Horizon $H$ is $2400$ days with possible visits every $15$, $30$ or $60$ days ($\delta=15$, $\overline\delta=60$). 

If treatment $\emptyset$ is applied, the patient may randomly jump from $m=0$ to any of the two disease states $m\in\{1,2\}$. In the sound state $m=0$, the marker level is contant. In the disease states ($m=1$ or $m=2$), the marker level follows an increasing flow, no other jump is possible.

If treatment $a$ is applied, the patient may only randomly jump from $m=0$ to the other disease state $m=2$. 
In the disease state $m=1$, the marker follows a decreasing flow and may randomly jump to the other disease state $m=2$.
In the disease state $m=2$, the marker follows an increasing flow and no other change of state is possible. Effects of treatment $b$ is similar: decrease of the marker in disease $m=2$, increase  in disease $m=1$.
The specific values of the local characteristics of the PDMP and proof of their verifying the assumptions of Section~\cref{sec:regularity} can be found in the supplementary materials, and the state graph for mode transitions is given in \cref{fig:states}.
%
\begin{figure}[htp]
\begin{center}
\begin{tikzpicture}
\node[draw,rounded corners=2pt] (mode0) at (0,0) {$0$};
\node[draw,rounded corners=2pt] (mode1) at (1.5,1) {$1$};
\node[draw,rounded corners=2pt] (mode2) at (1.5,-1) {$2$};

\draw[->,>=latex,color=blue!60] (mode2.north east) -- (mode1.south east) node[right=0.01cm,pos=0.5]{\small{$b$}};
\draw[->,>=latex,color=blue!60] (mode1.south west) -- (mode2.north west) node[left=0.01cm,pos=0.5]{\small{$a$}};
\draw[->,>=latex] (mode2) -- (mode0) node[right=0.01cm,pos=0.7]{\small{$b$}};
\draw[->,>=latex] (mode1) -- (mode0) node[right=0.01cm,pos=0.7]{\small{$a$}};
\draw[->,>=latex,color=red!60] (mode0.south west) -- (mode2.south west) node[left=0.01cm,pos=0.6]{\small{$\neq b$}};
\draw[->,>=latex,color=red!60] (mode0.north west) -- (mode1.north west) node[left=0.01cm,pos=0.6]{\small{$\neq a$}};
\end{tikzpicture}
\end{center}
\caption{\label{fig:states} {{\bf State graph for the changes of mode in the medical example.} Letters indicate under which treatments the jumps are possible. Black arrows indicate changes to remission phase, red arrows changes to \textit{standard} relapse, and blue arrows changes to \textit{therapeutic escape}}.}
\end{figure}
%
\subsection{Examples in biology}
\label{sec:bioexample}
%
{One may consider other kinds of therapeutic treatments, such as in \cite{pasin_controlling_2018} where the treatments correspond to interleukin injections in HIV-infected patients under antiretroviral treatment that must be suitably dosed and timed in order to maintain lymphocytes levels over a prescribed limit. In this example, the PDMP has no free modes and a controller-chosen corresponding to non-proliferating and proliferating lymphocytes, which dynamics can be influenced by punctual interleukin injections. Natural jumps correspond to the efficiency duration of an injection. Hence there is at most one jump between two injections. It is natural to suppose that the lymphocytes is observed through noisy proxies at discrete dates.}

{Other examples in biology may cover for instance gene expression or chemostats. The simplest PDMP models for gene expression include a (possibly controller triggered) promoter that can be on or off, and continuous levels of RNA and proteins, see e.g. \cite{herbach_inferring_2017}. If only the controller can turn the promoter on or off, there are no natural jumps. One may also consider a network of interacting genes with partial control over promoters and feedback loops. Chemostats are biorecactors used to cultivate some micro-organisms. They can be modeled by PDMPs, see e.g. \cite{fritsch2015modeling}, and are controlled by tuning the flow of fresh medium. Again, it is natural to consider hidden information in these frameworks.}
%
\subsection{Dynamic reliability}
\label{sec:fiabexample}
%
{The field of dynamic reliability is particularly suitable for controlled PDMP models, as explained in the introduction. The number of jumps in a given time interval is usually bounded for non repairable systems, for instance as soon as there is a finite number of components, each with a single nominal state and a finite number of successive possible degradation states up to total failure. It is also common not to observe the state of each component at all times.}

{The simplest examples in this field are pure jump Markov processes. It fits our assumptions as the flow and intensity are constant in every mode, and the state space is discrete hence all kernels are regular. See \cite{makis2003optimal} for such an example in the simpler case of optimal stopping. The next important subclass is when the mode corresponds to the states of the components, and the Euclidean variables to their age or use-time, see e.g. \cite{baysse2013maintenance,de2019dynamic} or \cite[Section 1.8.4]{dSDZ16} for a repair-workshop model with similar variables. In this case, the control $\ell$ may correspond to a change of settings.}

{A typical example of a more complex controlled PDMP in this field is the heated hold-up tank, a well know test case for dynamic reliability, see e.g. \cite{devooght97, MZ96, ZDDG08}. The system consists of a tank containing a fluid whose level is controlled by three components: two inlet pumps and one outlet valve. A thermal power source heats up the fluid. The failure rate of the components depends on the temperature, the position of the three components monitors the liquid level in the tank, and in turn, the liquid level determines the temperature. The pumps and valve can be either functioning or stuck, which corresponds to the free modes $m$, and turned on or off (if functioning), corresponding to the controller-chosen mode $\ell$. In any setting $\ell$, there are at most three jumps (corresponding to the failure of the pumps and valve). The flow and jump parameters have analytical expressions and satisfy our regularity assumptions.}
%
\section{Simulation study}
\label{sec:simulation}
%
We consider in this section the medical example of Section \ref{sec:medexample}. We briefly discuss the choice of the cost function and the grid construction, then state competing approaches and compare their performance to our approach, based on Monte Carlo simulations.
%
\subsection{Choice of the cost functions} 
\label{sec:cost}
%
The candidate policy depends on the cost functions, therefore the latter has to be carefully chosen. We denote $c_i(x,d)=C_V$ a fixed cost per visit that takes into account an emotional burden for the patient and health care expenses for the check-up. For the counterpart of the integral of the running cost, we choose en expression of the form
\begin{align*}
\tilde c(x,d,x')=\kappa |\zeta'-K_i|r+\beta r \1_{\{\ell\neq\emptyset\}}.
\end{align*} 
Parameter $\beta>0$ represents the cost of the treatment, $\kappa>0$ is a scale parameter and $ |\zeta'-K_i|r$ is a (crude) proxy of the integral of the process. Our cost function is then $c(x,d,x')=C_V+\tilde c(x,d,x')$. Calibrating those parameters is a difficult task that is not discussed here. 
%
\subsection{Construction of the grids}
\label{sec:grids}
%
Because decisions influence the dynamics of the process, grids cannot be constructed with techniques based on simulations, such as quantization. For instance, on a horizon of $2400$ days with possible visits every $15$, $30$ or $60$ days, this leads to approximately $10^{152}$ possible strategies. Grids therefore have to be chosen by expert knowledge, and transition probabilities computed accordingly. In the supplementary materials we discuss the construction of parcimonious grids that provide satisfying performance.
\subsection{Strategies in competition and performance criteria}
%
To evaluate the performance of our approach (OS in the comparison tables) we compare our work with several other strategies. 

The (unachievable) gold standard strategy is the \emph{See All}  (SA)
where decisions are taken  while observing the process $X_{t_n}$ by choosing the optimal treatment for the current mode. In this SA strategy, we do not allow choice of the next visit date. 

The \emph{filter} strategy corresponds to choosing the optimal treatment for the estimated mode: at each new observation $y_{n}$ the approximate filter $\bar\theta_{n}$ is computed based on the first discretization. Then the most probable mode $\hat m_{n}$ is obtained as
\begin{align*}
\hat m_{n}=\arg\max_{m=0:3}\sum_{j; \omega^j=(m,\cdot)}\bar\theta_{n}^j.
\end{align*}
Note that to set up this strategy, only the first discretization is needed. As this discretization is the least computationally extensive, it might be worth designing larger grids $\Omega$ to obtain better mode estimates. However, this strategy does not take into account knowledge on the dynamics of the process, hence we expect it to have higher costs. It is not able to select the next decision date either.

 The \emph{Standard} strategy is classically used in the hospitals and is based on a threshold $s_{rel}$ for relapse, and on a fixed duration for treatment $tlength$. While the observations remain below $s_{rel}$, the patient does not receive treatment and visits are scheduled every 2 months. When $s_{rel}$ is reached, the practitioner gives treatment $b$ (corresponding to the most frequent relapse type $2$) and the next visit is scheduled in $15$ days. If at the next visit the observed marker level has decreased, treatment $b$ is maintained with visits every $15$ for a total duration of $tlength$ days. Otherwise the practitioner changes treatment for treatment $a$ (with visits every $15$ days, and same total duration). Once treatment is stopped, visits are scheduled every two months and the strategy is repeated.

The \emph{fixed dates} optimal strategies are the strategies based on our discretization and dynamic programming approach where the choice of the next visit date are not allowed. We will investigate the fixed dates every $15$ and $60$ days (FD-15 and FD-60 in the comparison tables).

The performance criterium is the real cost of each strategy evaluated on the real process $X$ (averaged over $1000$ simulations). 
%
\subsection{Results} 
\label{Results}
%
First, note that 
all bits of codes are available online at 
\url{https://github.com/acleynen/PDMP-control}.
All results presented here are based on $1000$ simulations. The exact specifications 
can be found in the supplementary materials.

A first interesting result relates to the allocation of the computational burden. We compared, for a fixed number of elements of the second grid ($n_\Gamma$), the performance of our approach with different number of points in the first grid ($n_\Omega$), the smallest ones being included in the largest one. The performance of the approach drastically improves with a higher number of point in $\Omega$, leading to a very crude discretization of $\mathcal{P}(\Omega)$, compared to a less refined discretization of $\mathbb{X}$ but leading to a more refined discretization of $\mathcal{P}(\Omega)$. In some cases, adding elements to $\Gamma$ in fact worsened the results, as novel elements allowed to consider filters with less certainty in the values of the process, most often leading to worse decisions. This is illustrated in table 5 of the Supplementary Materials.

{We therefore compared our approach with competing strategies on the largest $\Omega$ grid, containing 133 elements in $\Omega$ and 144 elements in $\Gamma$. \cref{fig:traj} illustrates a controlled trajectory with treatments indicated by shape and disease indicated by color, for our approach. One can note that in most cases, the first time the projected filter has a  most likely element in a disease mode, treatment is not administered but the time lapse to the next visit is reduced. Only on the next visit will treatment be allocated. This explains the differences of performance between the different strategies given in \cref{tab:comp} and shows the added value of taking into account the knowledge on the dynamics of the process. It performs significantly better.}
\begin{figure}[tp]
  \begin{center}
  \includegraphics[width=0.9\linewidth]{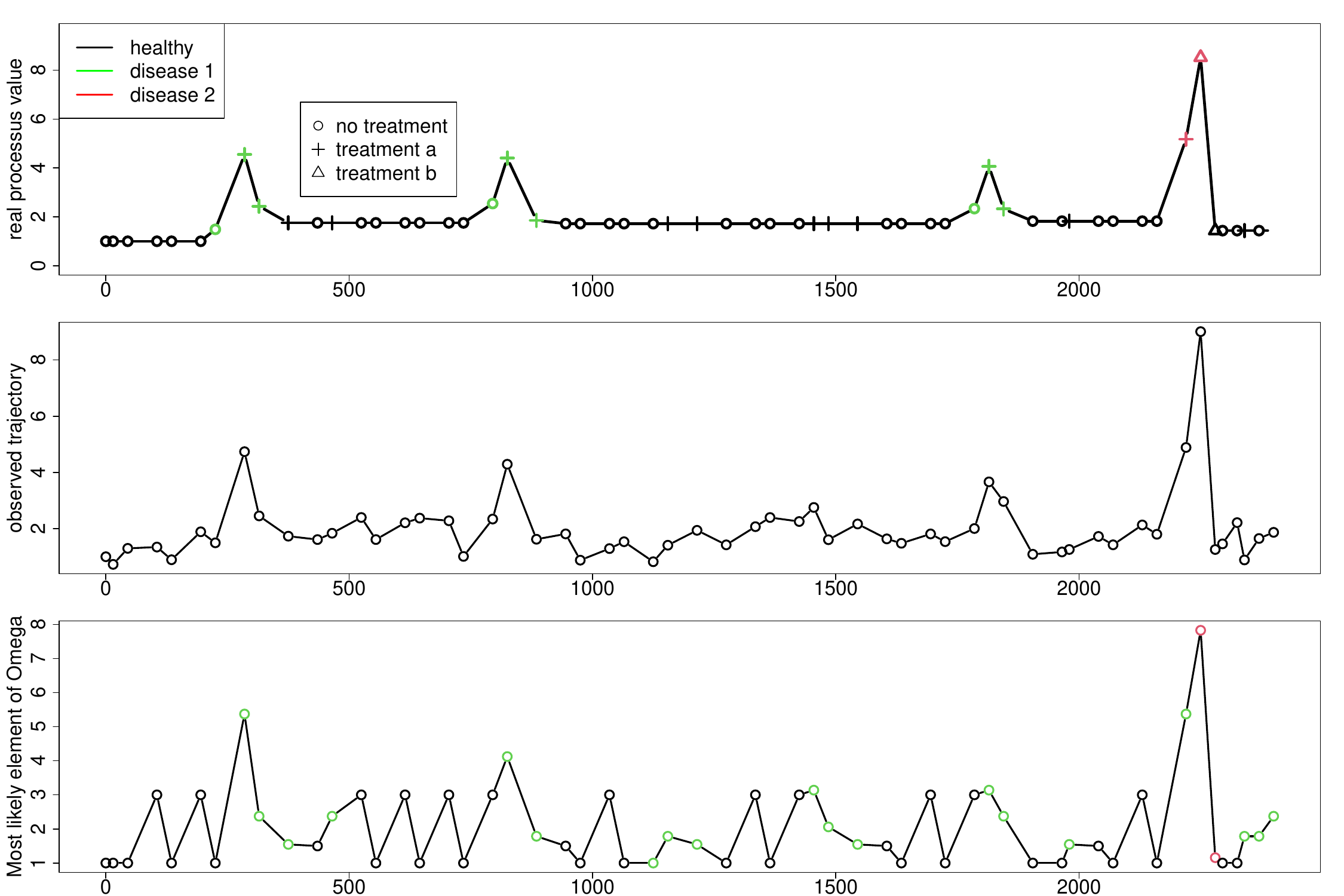}
    \caption{\label{fig:traj} \textbf{Example of a controlled (OS) simulated strategy.} The top panel illustrates the (hidden) values of the true process, with modes indicated by color, and allocated treatment indicated by point shape. The middle panel illustrates the observed trajectory (noisy version of the euclidean variable from the panel above). The bottom panel illustrates the most likely element of $\Omega$ estimated by the filter of the process.}
  \end{center}
  \end{figure}
\begin{table}[tp]
\footnotesize
\caption{\label{tab:comp} \textbf{Strategy comparisons.} Performance of each strategy measured in terms of Monte Carlo costs.}
\begin{center}
\begin{tabular}{r|r|r}
  \hline
 & Visits &  real cost (sd) \\
  \hline
 \hline
 & Choice & \textbf{203  (3.03)}   \\ 
  OS & 15 days & 291.36  (2.66)  \\ 
   & 60 days & 215.98  (3.2)  \\ 
   \hline
 & 15 jours & 321.28  (2.2)  \\ 
  Filter & 60 jours & 223.8  (4.03) \\ 
   \hline
 & 15 days & 203.15  (0.89)   \\ 
   See All   & 60 days & \textbf{113.72}  (1.79)     \\ 
   \hline
Standard &  & 301.09  (5.69)     \\ 
   \hline
\end{tabular}
\end{center}
\end{table}
%
\section{Conclusion} \label{sec:conclusion}

We have proposed a numerically feasible scheme to approximate the value function of an impulse control problem for a class of hidden PDMPs where control actions do influence the dynamic of the process. Approximations rely on discretization of the observed and unobserved state spaces, for which we propose strategies to explore the large-dimensional belief spaces that depend on the process characteristics while maintaining error bounds for this approximation explicitly depending on the parameters of the problem. Codes are made freely available on Github.

This provides a mathematical setting to study realistic processes for instance in a disease-control framework, which we illustrate on simulations. This is a promising start
for the  study of even more realistic disease-control frameworks. 

The important
open question concerns the optimality of our candidate
strategy. It cannot be directly linked to our various operators, but we are hopeful that further work will enable
us to prove theoretically that it is close to optimality.

\appendix
\section{Proofs of the main results}
\subsection{Stability of kernels $R'$ and $\bar{R}'$}
\subsubsection{Kernel $R'$ maps $\X'$ onto itself}
\label{app:R'X'}
%
Set $\xi\in\X'$ and $d\in\K'(\xi)$, and denote $\Xi$ a random variable with distribution $R'(\cdot|\xi,d)$.
If $d=\cd$ then $\Xi=\Delta\in\X'$.
Otherwise, by definition one has $\xi=(\theta,y,w)$ such that
$\theta(E_{m_0})\geq \Big(\frac{\underline f}{\overline f}\Big)^{\frac{w}{\delta}\vee 1}e^{-w\|\lambda\|}$, 
and $d=(\ell,r)\in \K'(\xi)$.
Then, one has 
 $\Xi=(\Theta,Y,W)$ and conditionally to $Y$, 
$
\Theta=\Psi(\theta,Y,w+r,d).
$
Now, the definition of the filter operator $\Psi$ and Assumption \ref{H1} yield
\begin{align*}
\Theta(E_{m_0})
&= \frac{\int_{E_M} \int_{E_M}f(Y-F(x'))\1_{E_{m_0}}(x')P(dx'|x,d)\theta(dx)}{\int_{E_M} \int_{E_M}f(Y-F(x'))P(dx'|x,d)\theta(dx)}\\
&\geq \frac{\underline{f}}{\overline{f}}\int_{E_{m_0}} P(E_{m_0}|x,d)\theta(dx).
\end{align*}
By definition of kernel $P$, if $x\in E_{m_0}$ and $d=(\ell,r)$, the probability to be in $E_{m_0}$ after a time $r$ and starting from $E_0$ at time $0$ is greater than the probability not to jump at all between times $0$ and $r$. Thus one has 
\begin{align*}
P(E_{m_0}|x,d)\geq e^{-\bs\Lambda^{\ell}(x,r)}\geq e^{-r\|\lambda\|},
\end{align*}
thanks to Assumption \ref{H-lambda-b}. Then one has
\begin{align*}
\Theta(E_{m_0})
&\geq \frac{\underline f e^{-r\|\lambda\|}\theta(E_{m_0})}{\overline f}\geq \Big(\frac{\underline f}{\overline f}\Big)^{\frac{w+r}{\delta}\vee 1}e^{-(w+r)\|\lambda\|},
\end{align*}
as $\frac{\underline f}{\overline f}\leq 1$ and $r\geq \delta$. 
%
\subsubsection{Kernel $\bar{R}'$ maps $\bar{\mathbb{X}}'$ onto itself}
\label{app:R'X'bar}
%
Set $\bar\xi\in\bar\X'$ and $d\in\K'(\xi)$, and denote $\Xi$ a random variable with distribution $\bar R'(\cdot|\bar\xi,d)$.
If $d=\cd$ then $\Xi=\Delta\in\bar\X'$.
Otherwise, by definition one has $\bar\xi=(\bar\theta,y,w)$ such that $\bar\theta(\Omega_{m_0})\geq \Big(\frac{\underline f}{\overline f}\Big)^{\frac{w}{\delta}\vee 1}e^{-w\|\lambda\|}$, 
and $d=(\ell,r)\in \K'(\xi)$.
In this case, one has $\Xi=(\Theta,Y,W)$ and conditionally to $Y$,  
$
\Theta=\Psibar(\bar\theta,Y,w+r,d),
$ 
and
\begin{align*}
\Theta(\Omega_{m_0})
  &= \frac{\sum_{\omega^i\in \Omega}\sum_{\omega^j\in \Omega_{m_0}} f(Y-F(\omega^j))\bar P(\omega^j|\omega^i,d)\thetabar(\omega^i)}{\sum_{\omega^i\in \Omega}\sum_{\omega^k\in \Omega}f(Y-F(\omega^k))\bar P(\omega^k|\omega^i,d)\thetabar(\omega^i)}\\
  &\geq \frac{\underline f}{\overline f}\sum_{\omega^i\in \Omega_{m_0}}\bar P(\Omega_{m_0}|\omega^i,d)\thetabar(\omega^i),
\end{align*}
by Assumption  \ref{H1} again. A similar argument as in the previous section directly yields
$\bar P(\Omega_{m_0}|\omega^i,d) = P(E_0|x,d)\geq e^{-r\|\lambda\|}$ 
as the Voronoi cells form a partition of $E_M$ that preserves the mode, and 
\begin{align*}
\Theta(\Omega_{m_0})
\geq \frac{\underline f e^{-r\|\lambda\|}\bar\theta(\Omega_{m_0})}{\overline f}
\geq \Big(\frac{\underline f}{\overline f}\Big)^{\frac{w+r}{\delta}\vee 1}e^{-(w+r)\|\lambda\|},
\end{align*}
hence the result.
%
\subsection{Error bounds for the first discretization}
\label{app:error1}
%
The proof of Theorem \ref{theorem1} is split into a number of intermediate propositions and lemmas stated in the sequel. Most of the proofs rely on adequate splitting of expressions into simpler terms, and exploitation of Lipschitz regularity. As they are sometimes long and partly straightforward, some details are omitted. We highlight the main ideas and point out  where we use our specific assumptions.
%
\subsubsection{Function spaces}
\label{app:function_spaces}
%
{We start by defining some function spaces in order to evaluate the regularity of our operators. 
Denote $BL(\X)$ the set of real-valued, bounded measurable functions $h$ on $\X$ that are Lipschitz continuous and set 
\begin{align*}
\|h\|=\sup_{\xi \in \X}\|h(\xi)\|,&\qquad
[h]=\sup_{m\in M}\sup_{\gamma\in \O}\sup_{\mathtt{x}\neq \mathtt{x'} \in E_m}\frac{|h(m,\mathtt{x},\gamma)-h(m,\mathtt{x}',\gamma)|}{\|\mathtt{x}-\mathtt{x}'\|}.
\end{align*}
%
Note also $BL_1(\X)=\{h\in BL(\X); \ \|h\|+[h]\leq 1\}$, and 
$BL(E)$ the restriction of $BL(\X)$ to functions defined on $E_M$.
Denote $BL(\X,\overline\delta)$ the set of real-valued, bounded measurable functions $w$ on $\X\times[0,\overline\delta]$ that are Lipschitz continuous and set 
\begin{align*}
\|w\|=\sup_{(\xi,t)\in \X\times\R}\|w(\xi,t)\|,\quad
[w]=\sup_{m\in M}\sup_{\gamma\in \O}\sup_{\mathtt{x}\neq \mathtt{x'} \in E_m\times\R}\frac{|w(m,\mathtt{x},\gamma,t)-w(m,\mathtt{x}',\gamma,t')|}{\|\mathtt{x}-\mathtt{x}'\|+|t-t'|}.
\end{align*}
%
Note also $BL(E,\overline\delta)$ the restriction of $BL(\X)$ to functions defined on $E_M\times\R$.
For $\theta$ and $\thetabar $ two probability measures in $\mathcal{P}(E_M)$, define the distance $d_E(\theta, \thetabar )$ by 
\begin{align*}
d_E(\theta, \thetabar )=\sup_{h \in BL_1(\X)}\sup_{\gamma\in\O}   \left|\int h(x,\gamma)\theta(dx) -\int h(x,\gamma) \thetabar  (dx)\right|.
\end{align*}
%
Let $BLP(\X')$ be the set of real-valued, bounded measurable functions $\varphi$ on $\X'$ that are Lipschitz continuous and set
\begin{align*}
\|\varphi\|=\sup_{\xi\in \X'}\|\varphi(\xi)\|,\quad
[\varphi]=\sup_{(\theta,\gamma)\neq (\theta',\gamma)\in \X'-\{\Delta\}}\frac{|\varphi(\theta,\gamma)-\varphi(\thetabar ,\gamma)|}{d_E(\theta,\thetabar )}.
\end{align*}
}
%
%
\subsubsection{Regularity of operator $P$}
\label{app:reg-P}
%
{The generic form of the transition kernels $P$ of the skeleton chains with time span $r$, in the special case where at most 2 jumps can occur before time $r$ is
\begin{align*}
Ph(x,d)
=&I^\ell[h](x,r)
+J^\ell\left[I^\ell[h]\right](x,r)
+J^\ell\left[J^\ell\left[H^\ell[h]\right]\right](x,r),
\end{align*}
%
where operator $H$ correspond to no jump, and no more possible jumps, $I$ to no jump in the case where more jumps can still occur, operator $J$ to one random jump. For any bounded measurable functions $h$ on $E_M$ and $w$ on $E_M\times [0,H]$, all $\ell\in L$, $x\in E_M$ and $0\leq t\leq \overline\delta$, these operators are defined as
%
\begin{align*}
H^\ell[h](x,t) & = h(\bs\Phi^\ell(x,t)),\\
I^\ell[h](x,t) & = h(\bs\Phi^\ell(x,t))e^{-\bs\Lambda^\ell(x,t)},\\
J^\ell[w](x,t) &= \int_{0}^{t}\int_{E_M} w(x',t-s) \bs Q^\ell (dx'|\bs\Phi^\ell(x,s))\bs \lambda^\ell(\bs\Phi^\ell(x,s))e^{-\bs\Lambda^\ell(x,s)}ds.
\end{align*}
%
If more jumps are allowed, the iterates are longer and the inner term is $I^\ell[h]$ for all iterates except the last one, and $H^\ell[h]$ for the last one: for at most $N_j$ jumps, one has
\begin{align*}
Ph(x,d)
=\sum_{k=0}^{N_j-1}(J^\ell)^k\left[I^\ell[h]\right](x,r) + 
(J^\ell)^{N_j}\left[H^\ell[h]\right](x,r).
\end{align*}
%
The sequel is detailed in the case of $N_j=2$ but can be readily extended to any other finite value. To ensure the regularity of $P$, one thus needs that the 3 operators above are regular enough, and that this property is preserved through iterations. 
\begin{lem}\label{lem-H}
Under assumption \ref{H-Phi}, for all $\ell\in L$, operator $H^\ell$ sends $BL(E)$ onto $BL(E,\overline\delta)$: for all function $h\in BL(E)$, for all  $x^1,x^2\in E_M$, $t^1,t^2\in[0,\overline\delta]$, one has
\begin{align*}
|H^\ell[h](x^1,t^1)|&\leq \|h\|,\\
|H^\ell[h](x^1,t^1)-H^\ell[h](x^2,t^2)|&\leq [h][\Phi](\|x^1-x^2\|+|t^1-t^2|).
\end{align*}  
Set $C_H=[\Phi]$.
\end{lem}
%
\begin{proof}
The first statement is obvious. Let us turn to the Lipschitz regularity.
For all function $h\in BL(E)$, one has
\begin{align*}
|H^\ell[h](x^1,t^1)-H^\ell[h](x^2,t^2)|
&= |h(\bs\Phi^\ell(x^1,t^1))-h(\bs\Phi^\ell(x^2,t^2))|\\
&\leq [h][\Phi](\|x^1-x^2\|+|t^1-t^2|),
\end{align*}  
by assumption \ref{H-Phi}.
\end{proof}
%
\begin{lem}\label{lem-I}
Under assumptions \ref{H-Phi} to \ref{H-lambda-b}, for all $\ell\in L$, operator $I^\ell$ sends $BL(E)$ onto $BL(E,\overline\delta)$: for all function $h\in BL(E)$, for all  $x^1,x^2\in E_M$, $t^1,t^2\in[0,\overline\delta]$, one has
\begin{align*}
|I^\ell[h](x^1,t^1)|&\leq \|h\|,\\
|I^\ell[h](x^1,t^1)-I^\ell[h](x^2,t^2)|&\leq ([h][\Phi]+\|h\|\overline\delta[\lambda])\|x^1-x^2\|+([h][\Phi]+\|h\|\|\lambda\|)|t^1-t^2|.
\end{align*}  
Set $C_I=[\Phi]+\overline\delta[\lambda]+[\Phi]+\|\lambda\|$.
\end{lem}
%
\begin{proof}
The first statement is obvious. Let us turn to the Lipschitz regularity.
For all function $h\in BL(E)$, one has
\begin{align*}
\lefteqn{|I^\ell[h](x^1,t^1)-I^\ell[h](x^2,t^2)|}\\
&= |h(\bs\Phi^\ell(x^1,t^1))e^{-\bs\Lambda^\ell(x^1,t^1)}-h(\bs\Phi^\ell(x^2,t^2))e^{-\bs\Lambda^\ell(x^2,t^2)}|\\
&\leq |H^\ell[h](x^1,t^1)-H^\ell[h](x^2,t^2)|e^{-\bs\Lambda^\ell(x^1,t^1)}+|H^\ell[h](x^2,t^2)|\left|e^{-\bs\Lambda^\ell(x^1,t^1)}-e^{-\bs\Lambda^\ell(x^2,t^2)}\right|\\
&\leq [h][\Phi](\|x^1-x^2\|+|t^1-t^2|)+\|h\|(\overline\delta[\lambda]\|x^1-x^2\|+\|\lambda\||t^1-t^2|),
\end{align*}  
by Lemma \ref{lem-H} and assumptions \ref{H-Phi}, \ref{H-lambda-lip}, \ref{H-lambda-b}.
\end{proof}
%
\begin{lem}\label{lem-J}
Under assumptions \ref{H-Phi} to \ref{H-Q}, for all $\ell\in L$, operator $J^\ell$ sends $BL(E,\overline\delta)$ onto itself: for all function $w\in BL(E,\overline\delta)$, for all  $x^1,x^2\in E_M$, $t^1,t^2\in[0,\overline\delta]$, one has
\begin{align*}
|J^\ell[w](x^1,t^1)|&\leq \|w\|,\\
|J^\ell[w](x^1,t^1)-J^\ell[w](x^2,t^2)|&\leq ([Q][w][\Phi]+\|w\|\overline\delta([\lambda]+ \overline\delta\|\lambda\|))\|x^1-x^2\|+2[w]|t^1-t^2|.
\end{align*}  
Set $C_J=[Q][\Phi]+\overline\delta([\lambda]+ \overline\delta\|\lambda\|)+2$.
\end{lem}
%
\begin{proof}
The first statement is obvious. Let us turn to the Lipschitz regularity. Suppose without loss of generality that $t^1\leq t^2$. 
For all function $w\in BL(E,\overline\delta)$, one has
\begin{align*}
\lefteqn{|J^\ell[w](x^1,t^1)-J^\ell[w](x^2,t^2)|}\\
&= \Big| \int_{0}^{t^1}\int_{E_M} w(x',t^1-s) \bs Q^\ell (dx'|\bs\Phi^\ell(x^1,s))\bs \lambda^\ell(\bs\Phi^\ell(x^1,s))e^{-\bs\Lambda^\ell(x^1,s)}ds\\
&\quad - \int_{0}^{t^2}\int_{E_M} w(x',t^2-s) \bs Q^\ell (dx'|\bs\Phi^\ell(x^2,s))\bs \lambda^\ell(\bs\Phi^\ell(x^2,s))e^{-\bs\Lambda^\ell(x^2,s)}ds\Big|\\
&\leq \Big| \int_{0}^{t^1} \Big| \int_{E_M} w(x',t^1-s) \bs Q^\ell (dx'|\bs\Phi^\ell(x^1,s))\bs \lambda^\ell(\bs\Phi^\ell(x^1,s))e^{-\bs\Lambda^\ell(x^1,s)}\\
&\quad - \int_{E_M} w(x',t^2-s) \bs Q^\ell (dx'|\bs\Phi^\ell(x^2,s))\bs \lambda^\ell(\bs\Phi^\ell(x^2,s))e^{-\bs\Lambda^\ell(x^2,s)} \Big| ds\\
&\quad + \int_{t^1}^{t^2}\int_{E_M} |w(x',t^2-s)| \bs Q^\ell (dx'|\bs\Phi^\ell(x^2,s))\bs \lambda^\ell(\bs\Phi^\ell(x^2,s))e^{-\bs\Lambda^\ell(x^2,s)}ds\\
&\leq A+ B.
\end{align*}  
The second term $B$ is readily bounded by $|t^1-t^2|\|w\|$. Let us turn to the first term. One has
\begin{align*}
A 
&\leq \int_{0}^{t^1} \Big| \int_{E_M} w(x',t^1-s) \bs Q^\ell (dx'|\bs\Phi^\ell(x^1,s))-\int_{E_M} w(x',t^2-s) \bs Q^\ell (dx'|\bs\Phi^\ell(x^2,s))\Big|\\
&\quad \times\bs \lambda^\ell(\bs\Phi^\ell(x^1,s))e^{-\bs\Lambda^\ell(x^1,s)}ds\\
&\quad +   \int_{0}^{t^1} \int_{E_M} |w(x',t^2-s)| \bs Q^\ell (dx'|\bs\Phi^\ell(x^2,s))\Big|\bs \lambda^\ell(\bs\Phi^\ell(x^1,s))e^{-\bs\Lambda^\ell(x^1,s)} - \bs \lambda^\ell(\bs\Phi^\ell(x^2,s))e^{-\bs\Lambda^\ell(x^2,s)} \Big| ds\\
&\leq A_1+A_2.
\end{align*}  
For term $A_1$, one has
\begin{align*}
A_1
&\leq \int_{0}^{t^1} \int_{E_M} |w(x',t^1-s)-w(x',t^2-s)| \bs Q^\ell (dx'|\bs\Phi^\ell(x^1,s))\bs \lambda^\ell(\bs\Phi^\ell(x^1,s))e^{-\bs\Lambda^\ell(x^1,s)}ds\\
&\quad + \int_{0}^{t^1} \Big| \int_{E_M} w(x',t^2-s) \bs Q^\ell (dx'|\bs\Phi^\ell(x^1,s))-\int_{E_M} w(x',t^2-s) \bs Q^\ell (dx'|\bs\Phi^\ell(x^2,s))\Big|\\
&\quad \times\bs \lambda^\ell(\bs\Phi^\ell(x^1,s))e^{-\bs\Lambda^\ell(x^1,s)}ds\\
&\leq A_{11}+ A_{12}.
\end{align*}  
On the one hand, one has
\begin{align*}
A_{11}
&\leq [w]|t^1-t^2|\int_{0}^{t^1} \bs \lambda^\ell(\bs\Phi^\ell(x^1,s))e^{-\bs\Lambda^\ell(x^1,s)}ds
\leq [w]|t^1-t^2|.
\end{align*}  
On the other hand, 
for fixed $t^2,s$, the mapping $h_{t^2,s}:x\mapsto w(x,t^2-s)$ is clearly in $BL(E)$, with $\|h_{t^2,s}\|\leq \|w\|$ and $[h_{t^2,s}]\leq [w]$. Hence, assumption \ref{H-Q} yields
\begin{align*}
A_{12}
&= \int_{0}^{t^1} \Big| \bs Q^\ell h_{t^2,s}(\bs\Phi^\ell(x^1,s))-\bs Q^\ell h_{t^2,s}(\bs\Phi^\ell(x^2,s))\Big|\bs \lambda^\ell(\bs\Phi^\ell(x^1,s))e^{-\bs\Lambda^\ell(x^1,s)}ds\\
&\leq [Q][w]\int_{0}^{t^1} \Big|\bs\Phi^\ell(x^1,s)-\bs\Phi^\ell(x^2,s)\Big|\bs \lambda^\ell(\bs\Phi^\ell(x^1,s))e^{-\bs\Lambda^\ell(x^1,s)}ds\\
&\leq [Q][w][\Phi]\|x^1-x^2\|\int_{0}^{t^1} \bs \lambda^\ell(\bs\Phi^\ell(x^1,s))e^{-\bs\Lambda^\ell(x^1,s)}ds\ \leq [Q][w][\Phi]\|x^1-x^2\|,
\end{align*}  
by assumption \ref{H-Phi}. Let us now consider term $A_2$. One has
\begin{align*}
A_2
&= \int_{0}^{t^1} \int_{E_M} |w(x',t^2-s)| \bs Q^\ell (dx'|\bs\Phi^\ell(x^2,s))\Big|\bs \lambda^\ell(\bs\Phi^\ell(x^1,s))e^{-\bs\Lambda^\ell(x^1,s)} - \bs \lambda^\ell(\bs\Phi^\ell(x^2,s))e^{-\bs\Lambda^\ell(x^2,s)} \Big| ds\\
&\leq \|w\| \int_{0}^{t^1}\Big|\bs \lambda^\ell(\bs\Phi^\ell(x^1,s))e^{-\bs\Lambda^\ell(x^1,s)} - \bs \lambda^\ell(\bs\Phi^\ell(x^2,s))e^{-\bs\Lambda^\ell(x^2,s)} \Big| ds\\
&\leq \|w\|\overline\delta[\lambda]\|x^1-x^2\|+\|w\|\overline\delta^2\|\lambda\|\|x^1-x^2\|,
\end{align*}  
by assumptions \ref{H-Phi}, \ref{H-lambda-lip} and \ref{H-lambda-b}. The result follows.
\end{proof}
%
\begin{theo}\label{th-P-noboundary}
Under assumptions \ref{H-Phi} to \ref{H-jumps}, operator $P$ sends $BL(E)$ onto itself. More specifically, for all $d=(\ell,r)\in L\times\T$, for all function $h\in BL_1(E)$, for all  $x^1,x^2\in E_M$, one has
\begin{align*}
|P h(x^1,d)|&\leq \|h\|,\\
|P h(x^1,d))-Ph(x^2,d)|& \leq C_P\|x^1-x^2\|,
\end{align*}  
where $C_P$ is a constant depending only on the parameters of the problem.
\end{theo}
%
\begin{proof}
The first statement is obvious as $P$ is a Markov kernel. 
For all $x^1,x^2$ in $E_M$, one has
\begin{align*}
|P h(x^1,d))-Ph(x^2,d)|
&\leq |I^\ell[h](x^1,r)-I^\ell[h](x^2,r)|\\
&\quad +|J^\ell\left[I^\ell[h]\right](x^1,r)-J^\ell\left[I^\ell[h]\right](x^2,r)|\\
&\quad +|J^\ell\left[J^\ell\left[H^\ell[h]\right]\right](x^1,r)-J^\ell\left[J^\ell\left[H^\ell[h]\right]\right](x^1,r)(x^2,r)|\\
&\leq A+B+C.
\end{align*}
Lemma \ref{lem-I} applied to $h\in BL(E)$, yields
\begin{align*}
A
&= |I^\ell[h](x^1,r)-I^\ell[h](x^2,r)|
\leq ([h][\Phi]+\|h\|\overline\delta[\lambda])\|x^1-x^2\|.
\end{align*}
The same lemma yields that the mapping $w:(x,t)\mapsto I^\ell[h](x,t)$ is in $BL(E,\overline\delta)$, with $\|w\|\leq \|h\|$, and $[w]\leq C_I(\|h\|+[h])$. Then in turn, Lemma \ref{lem-J} applied to $w$ yields
\begin{align*}
B
&= |J^\ell\left[I^\ell[h]\right](x^1,r)-J^\ell\left[I^\ell[h]\right](x^2,r)|\ =  |J^\ell\left[w\right](x^1,r)-J^\ell\left[w\right](x^2,r)|\\
&\leq ([Q][w][\Phi]+\|w\|\overline\delta([\lambda]+ \overline\delta\|\lambda\|))\|x^1-x^2\|.
\end{align*}
The combination of Lemma \ref{lem-H} and \ref{lem-J} also yield that the mapping $w':(x,t)\mapsto J^\ell\left[H^\ell[h]\right](x,t)$ is still in $BL(E,\overline\delta)$, with $\|w'\|\leq \|H^\ell[h]\| \leq \|h\|$ and 
\begin{align*}
[w']&\leq C_J(\|H^\ell[h]\|+[H^\ell[h]])\leq C_J(\|h\|+C_H(\|h\|+[h])=(C_J+C_H)\|h\|+C_H[h].
\end{align*}
Lemma \ref{lem-J} for $w'$ now yields
\begin{align*}
C
&= |J^\ell\left[J^\ell\left[H^\ell[h]\right]\right](x^1,r)-J^\ell\left[J^\ell\left[H^\ell[h]\right]\right](x^1,r)(x^2,r)|\\
&= |J^\ell\left[w'\right](x^1,r)-J^\ell\left[w']\right](x^1,r)(x^2,r)|\\
&\leq  ([Q][w'][\Phi]+\|w'\|\overline\delta([\lambda]+ \overline\delta\|\lambda\|))\|x^1-x^2\|.
\end{align*}
hence the result with 
$C_P\leq \|h\|([Q][\Phi](C_H+C_I+C_J)+3\overline\delta[\lambda]+2\overline\delta^2\|\lambda\|))
+[h][\Phi](1+[Q](C_I+C_H))$.
\end{proof}
%
Note that further iterations of Lemmas \ref{lem-J} also yield that operator $P$ sends $BL(E)$ onto itself for any bounded number of allowed jumps in the interval $[0,r]$. 
Note also that Theorem \ref{th-P-noboundary} is still valid on $BL(\X)$.
}
%
\subsubsection{Regularity of operator $R'$}
\label{app:reg-R'}
%
We next need to investigate the regularity of the filter operator. The proof is based on the explicit form of the filter operator. The lower bounds on $\theta(E_{m_0})$ are required to bound the terms in the denominators.
%
\begin{lem}\label{Psi-theta-thetabar}
Under Assumptions \ref{H-Phi} to \ref{H-jumps}, there exist some positive constant $C_\Psi$ such that for all 
$w'\in[0,H]$, $d=(\ell,r) \in L\times\T$,  and $\theta,\thetabar \in \mathcal{P}(E_{M})$ satisfying $\theta(E_{m_0})>0$, $\thetabar(E_{m_0})>0$, we have
\begin{align*}
\int_I d_E(\Psi(\theta,y',w',d),\Psi(\thetabar,y',w',d))dy' \leq \frac{C_\Psi}{\theta(E_{m_0})\thetabar(E_{m_0})} d_E(\theta,\thetabar).
\end{align*}
\end{lem}
%
\begin{proof} Let  $ g \in BL_1(\X)$ and $y'\in I$, $w'\in[0,H]$, and $\gamma\in\O$, we have
\begin{align*}
\lefteqn{d_E(\Psi(\theta,y',w',d),\Psi(\thetabar,y',w',d))}\\
&=\int_{E_M} g(x,\gamma)\Psi(\theta,y',w',d)(dx)-  \int_{E_M} g(x,\gamma)\Psi(\thetabar,y',w',d)(dx)\\
&= \left(\int_{E_{M}^2}f(y'-F(x'))g(x',\gamma) P(dx'|x,d)\theta(dx)\right. \\
&\quad \left.-\int_{E_{M}^2} f(y'-F(x'))g(x',\gamma) P(dx'|x,d)\thetabar(dx)\right)
\frac{1}{\int_{E_{M}^2} f(y'-F(x'))P(dx'|x,d)\theta(dx)}\\
&\quad +\int_{E_{M}^2} f(y'-F(x'))g(x',\gamma) P(dx'|x,d)\theta(dx)\\
&\quad\times \left( \frac{1}{\int_{E_{M}^2} f(y'-F(x'))P(dx'|x,d)\theta(dx)} 
-  \frac{1}{\int_{E_{M}^2}   f(y'-F(x'))P(dx'|x,d)\thetabar(dx)}\right)\\
&=\frac{A}{B}+C.
\end{align*}
First, as $g$ is in $BL_1(\X)$, then 
$\varphi : (x,\gamma)\mapsto  \int_{E_M} f(y'-F(x')) g(x',\gamma)P(dx'|x,d)$ is in $BL(\X)$ with $\|\varphi\|\leq \overline f$ and $[\varphi]\leq \big(\overline f+L_Y\big)C_P$ thanks to Theorem \ref{th-P-noboundary}, so that one has
\begin{align*}
|A|&\leq \big(\overline f(1+C_P)+L_YC_P \big)d_E(\theta,\thetabar).
\end{align*}
On the other hand, one has
\begin{align*}
|B|&\geq \underline{f}\int_{E_{M}}P(E_{M}|x,d)\theta(dx)\geq \underline{f}\int_{E_{m_0}}P(E_{m_0}|x,d)\theta(dx).
\end{align*}
As seen in Section \ref{app:R'X'}, if $x\in E_0$ and $d=(\ell,r)$, one has 
\begin{align*}
P(E_0|x,d)=e^{-\bs\Lambda^{\ell}(x,r)}\geq e^{-r\|\lambda\|}\geq e^{-\overline\delta\|\lambda\|},
\end{align*}
hence, one has 
$
|B|\geq \underline{f}e^{-\overline\delta\|\lambda\|}\theta(E_{m_0}).
$ 
%
Similarly, the last term can be bounded by
\begin{align*}
|C|&\leq \overline{f}\frac{\overline f(1+C_P)+L_YC_P }{\underline{f}^2e^{-2\overline\delta\|\lambda\|}\theta(E_{m_0})\thetabar(E_{m_0})}d_E(\theta,\thetabar).
\end{align*}
As $I$ is a bounded interval, the result follows.
\end{proof}
%
Now we can turn to the regularity of operator $R'$.
%
\begin{prop}\label{Rf-BL}
Under Assumptions \ref{H-Phi} to \ref{H1}, there exists a positive constant $C_{R'}$ such that for all $g \in BLP_1(\X')$, $\gamma=(y,w) \in \O$, $d \in \K'(\gamma)$ and for all $\theta, \thetabar $ in $\mathcal P(E_M)$ {such that $(\theta,\gamma)\in\X'$, $(\thetabar,\gamma)\in\X'$,
one has
\begin{align*}
\big|R'g(\theta,y,d)-R'g(\thetabar,y,d)\big|
& \leq C_{R'}d_E(\theta,\thetabar ).
\end{align*}
}
In particular, operator $R'$ maps $BLP(\X')$ onto itself and for $g\in BLP(\X')$, one has $\|R'g\|\leq \|g\|$ and $[R'g]\leq (\|g\|+[g])C_{R'}$.
\end{prop}
\begin{proof}
If $d=\cd$, then $R'g(\theta,\gamma,\cd)=R'g(\thetabar ,\gamma,\cd)=g(\Delta)$, hence $|R'g(\theta,\gamma,\cd)-R'g(\thetabar ,\gamma,\cd)|=0$.
If $d=(\ell,r)\neq \cd$, then by assumption one has $\theta(E_{m_0})>0$ and $\thetabar (E_{m_0})>0$. In addition, one has  $w<N\delta$, such that 
\begin{align*}
R'g(\theta,\gamma,d)
=\int_{I}\int_{E_{M}^2} g\big(\Psi(\theta,y',w+r,d),y',w+r\big)f(y'-F(x'))P(dx'|x,d)\theta(dx)dy'.
\end{align*}
Set
\begin{align*}
h_\theta(x,\gamma)= \int_{I}g\big(\Psi(\theta,y',w+r,d),y',w+r\big)f(y'-F(x))dy',
\end{align*}
if $x\in E_{M}$ and $\gamma=(y,w+r)$ and $0$ otherwise and ${g}_{\theta,d}(x,\gamma)= Ph_\theta(x,\gamma,d)$. From Assumption \ref{H1}, $h_\theta$ is in $BL(\X)$ with $\|h_\theta\|\leq B_f$ and $[h_\theta]\leq L_f$. Theorem \ref{th-P-noboundary} thus yields that ${g}_{\theta,d}$ in still in $BL(\X)$ with $\|{g}_{\theta,d} \| \leq B_f$ and $[{g}_{\theta,d}]\leq (B_f+L_f)C_P$.
Therefore one has
\begin{align*}
\lefteqn{|R'g(\theta,\gamma,d)-R'g(\thetabar,\gamma,d)|}\\
&\leq  \int_{E_{M}}  \left|{g}_{\theta,d}(x,\gamma)-{g}_{\thetabar ,d}(x,\gamma)\right|\thetabar (dx) +  \left|\int_{E_{M}}  {g}_{\theta,d}(x,\gamma)\theta(dx)-\int_{E_{M}}  {g}_{\theta,d}(x,y) \thetabar (dx)\right| \\
&\leq  \int_{E_{M}} \int_I d_E(\Psi(\theta,y',w+r,d), \Psi(\thetabar ,y',w+r,d))dy'\thetabar (dx)+\left([{g}_{\theta,d}]+\|{g}_{\theta,d}\|\right) d_E(\theta,\thetabar ),
\end{align*}
as $g \in BLP_1(\X')$.  
Note in particular the assumptions on $\theta$ and $\thetabar$ guarantee that the assumptions of Lemma \ref{Psi-theta-thetabar} are satisfied. Hence we conclude using Lemma \ref{Psi-theta-thetabar} and the minoration of $\theta(E_{m_0})$ and $\thetabar(E_{m_0})$ from the definition of $\X'$.
\end{proof}
%
\subsubsection{Projection error for operators $P$ and $\bar P$}
\label{app:projP}
%
For $\omega=(m,\mathtt{x})\in \Omega$, we denote $\omega_1=m$, $\omega_2=\mathtt{x}$. 
For any function $h\in BL(\X)$, $\xi=(x,\gamma)\in\Omega\times\O$ and $d\in\K(\xi)$, denote
\begin{align*}
\bar Ph(\xi,d)= \int_{E_M} h(x',\gamma)\bar P(dx'|x,d).
\end{align*}
%
\begin{prop}\label{prop:P-Pbar}
For all $i \in \{1,\dots,{n_{\Omega}}\}$, $\gamma\in\O$, $d \in L\times\T$ and all $g\in BL_1(\X)$, we have 
\begin{align*}\label{eq:P-Pbar}
\Big|Pg(\omega^i,\gamma,d)-\bar{P}g(\omega^i,\gamma,d)\Big|\leq \sup_{j \in \{1,\dots,{n_{\Omega}}\}} \mathcal{D}_j.
\end{align*}
\end{prop}

\begin{proof}
For $x\in E$, denote by $C_x$ the Voronoi cell of $x$.
For $g\in BL_1(\X)$, $i \in \{1,\dots,{n_{\Omega}}\}$, $\gamma=(y,w)\in\O$, $d \in L\times\T$, one has
\begin{align*}
\Big|Pg(\omega^i,\gamma,d)-\bar{P}g(\omega^i,\gamma,d)\Big|
&\leq \Big|  \int_{E_M} g(x',\gamma') P(dx' | \omega^i, d)  -  \sum_{j=1}^{{n_{\Omega}}} g(\omega^{j},\gamma') P(C_j|  \omega^i, d)\Big| \\
& \leq   \sum_{j=1}^{{n_{\Omega}}}  \int_{C_{j}} \big| g(x',\gamma') -g(\omega^j,\gamma')\big| P(dx' | \omega^i, d)  \\
& \leq   \sum_{j=1}^{n_{\Omega}}  \int_{C_{j}} | x' -\omega^j |P(dx' | \omega^i, d)\ \leq  \sup_{j\in\{1,\dots,{n_{\Omega}}\}} \mathcal{D}_j ,
\end{align*}
as $g$ is Lipschitz-continuous in its first variable.
\end{proof}
%
\subsubsection{Projection error for operators $R'$ and $\bar R'$}
\label{app:projR}
%
We first need to evaluate the error between $\Psi$ and $\Psibar$.
\begin{lem}\label{dPsiPsibar}
Under Assumptions \ref{H-Phi} to \ref{H1}, there exist some positive constant $\bar C_\Psi$ such that for all $\gamma'=(y',w')\in\O$, $d=(\ell,r) \in L\times\T$ and  all {$\thetabar \in \mathcal{P}(\Omega)$ satisfying $\thetabar(\Omega_{m_0})>0$, we have
\begin{align*}
\int_I d_E(\Psi(\thetabar,\gamma',d),\Psibar(\thetabar,\gamma',d))dy' \leq \frac{\bar C_\Psi}{\thetabar(\Omega_{m_0})} \sup_{j \in \{1,\dots,\ell\}} \mathcal{D}_j.
\end{align*}
}
\end{lem}
%
\begin{proof}
Set $\gamma\in\O$, $\gamma'=(y',w')\in\O$, $d=(\ell,r) \in L\times\T$, $\thetabar \in \mathcal{P}(\Omega)$ satisfying  $\thetabar(\Omega_{m_0})>0$, and $g \in BL_1(\X)$. 
%
We have
\begin{align*}
\lefteqn{ \int_{E_M} g(x,\gamma)\Psi(\thetabar,y',w',d)(dx)- \int_{E_M} g(x,\gamma)\Psibar(\thetabar,y',w',d)(dx)}\\
&=\frac{\sum\limits_{\omega^i\in \Omega } \int_{E_{M}}f(y'-F(x'))g(x',\gamma) P(dx'|\omega^i,d)\thetabar(\omega^i)}{\sum\limits_{\omega^i\in \Omega }\int_{E_{M}}   f(y'-F(x'))P(dx'|\omega^i,d)\thetabar(\omega^i)} \\
&\cutline -\frac{\sum\limits_{\omega^i\in \Omega }\sum\limits_{\omega^j\in \Omega } f(y'-F(\omega^j))  g(\omega^j,\gamma)\bar P(\omega^j|\omega^i,d)\thetabar(\omega^i)}{\sum\limits_{\omega^i\in \Omega }\sum\limits_{\omega^k\in \Omega }{ f(y'-F(\omega^k))}\bar P(\omega^k|\omega^i,d)\thetabar(\omega^i)}.
\end{align*}
Using a similar splitting and similar arguments as in the proof of Lemma \ref{Psi-theta-thetabar} together with
Proposition \ref{prop:P-Pbar} yield the expected result, as  
$\varphi : (x,\gamma)\mapsto  \int_{E_M} f(y'-F(x')) g(x',\gamma)P(dx'|x,d)$ is still in $BL(\X)$ with $\|\varphi\|\leq \overline f$ and $[\varphi]\leq \big(\overline f+L_Y\big)C_P$ thanks to Theorem \ref{th-P-noboundary}.
\end{proof}
%
\begin{prop}\label{prop:R-Rbar}
Under Assumptions \ref{H-Phi} to \ref{H1}, there exist a positive constant ${\bar C_{R'}}$ such that for all $g \in BLP_1(E)$, {$(\theta,\gamma) \in \bar \X'$ and $d\in\K'(\gamma)$, one has
\begin{align*}
\left|R'g(\thetabar,\gamma,d)-\bar R'g(\thetabar,\gamma,d)\right| \leq {\bar C_{R'}} \sup_{j \in \{1,\dots,\ell\}} \mathcal{D}_j.
\end{align*}
}
\end{prop}
%
\begin{proof}
This is a direct consequence of Propositions \ref{prop:P-Pbar} and \ref{dPsiPsibar} with $\bar C_{R'}= L_f+2B_f+L_f\bar C_\Psi$, as $\thetabar(\Omega_{m_0})$ is uniformly bounded from below.
\end{proof}
%
\subsubsection{Regularity and approximation error for the cost functions}
\label{app:proof-cC}
%
The last preliminary result we need in order to prove Theorem \ref{theorem1} is to ensure that the cost functions $c'$, $C'$ belong to the appropriate function spaces and the error between $c'$ and $\bar c'$ is controlled. 
%
\begin{lem}\label{lem-C'}
Under Assumption \ref{Hc}, the cost function $C'$ is in $BLP(\X')$.
\end{lem}
%
\begin{proof}
Under Assumption \ref{Hc}, $C$ is clearly in $BL(E)$ with $\|C\|\leq B_C$ and $[C]\leq L_C$. 
Thus, $C'$ is still bounded by $B_C$ and for all $(\theta,\gamma)$ and $(\thetabar ,\gamma)$ in $\X'$ one has
\begin{align*}
|C'(\theta,\gamma)-C'(\thetabar,\gamma)|
&=\Big| \int_{E_M} C(x)\theta(dx)- \int_{E_M} C(x)\thetabar(dx)\Big|\leq (B_C+L_C)d_E(\theta,\thetabar).
\end{align*}
Thus $C'$ is in $BLP(\X)$ with $\|C'\|\leq B_C$ and $[C']\leq B_C+L_C$.
\end{proof}
%
\begin{lem}\label{lem-c'}
Under Assumptions \ref{H-Phi} to \ref{H-jumps} and \ref{Hc}, for all $d\in \A$ the function $c'_d:\xi\mapsto c'(\xi,d)$ is in $BLP(\X')$.
\end{lem}
%
\begin{proof}
If $d=\cd$, then $c'=C'$ and the result is true. Otherwise, under Assumption \ref{Hc}, $c'_d$ is bounded by $B_c$, and for fixed $x$, application $c_{x,d}:x'\mapsto c(x,d,x')$ is in $BL(E)$ with $\|c_{x,d}\|\leq B_c$ and $[c_{x,d}]\leq L_c$.
Let us study the application $Pc : x\mapsto  \int_{E_M} c(x,d,x')P(dx'|x,d)$ for fixed $d\in L\times\T$. It is still bounded by $B_c$, and for $x^1,x^2$ in $E$, one has
\begin{align*}
|Pc(x^1)-Pc(x^2)|
&\leq  \int_{E_M} |c(x^1,d,x')-c(x^2,d,x')|P(dx'|x^1,d) \\
& \quad +\Big| \int_{E_M} c(x^2,d,x')P(dx'|x^1,d)- \int_{E_M} c(x^2,d,x')P(dx'|x^2,d)\Big|\\
&\leq L_c\|x^1-x^2\|+C_P(B_c+L_c)\|x^1-x^2\|
\end{align*}
by applying Theorem \ref{th-P-noboundary} to $c_{x^2,d}$. Hence application $Pc$ is still in  $BL(E)$ hence in $BL(\X)$. 
Then, for all $(\theta,\gamma)$ and $(\thetabar ,\gamma)$ in $\X'$ one obtains
\begin{align*}
|c'_d(\theta,\gamma)-c'_d(\thetabar,\gamma)|
&=\Big|\int_{E_M^2}\!\!\! c(x,d,x')P(dx'|x,d)\theta(dx)-\int_{E_M^2}\!\!\!  c(x,d,x')P(dx'|x,d)\thetabar(dx)\Big|\\
&\leq (B_c+L_c)(C_P+1)d_E(\theta,\thetabar).
\end{align*}
Hence the result.
\end{proof}
%
\begin{lem}\label{lem-c'barc'}
Under Assumption \ref{Hc}, for all $\xi\in\bar \X'$, $d\in\K'(\xi)$, one has
\begin{align*}
|c'(\xi,d)-\bar c'(\xi,d)|\leq (B_c+L_c) \sup_{j \in \{1,\dots,\ell\}} \mathcal{D}_j.
\end{align*}
\end{lem}
%
\begin{proof}
The result the follows directly form Proposition \ref{prop:P-Pbar} and the facts that $x\mapsto c(\omega^i,d,x)$ is in $BL(E)$ and $c'(\cdot,\cd)=\bar c'(\cdot,\cd)$.
\end{proof}
\subsubsection{Proof of Theorem \ref{theorem1}}
\label{app:proof-th1}
%
We establish the result by (backward) induction on $n$, with the additional statement that $v'_n$ is in $BLP(\X')$ for all $n$.

$\bullet$ For $n=N$, $v'_N=C'$ is in $BLP(\X')$ by Lemma \ref{lem-C'}. In addition, set  $\xi$ in $\bar\X'$. Then by definition one has
\begin{align*}
|v'_N(\xi)-\bar v'_N(\xi)|=|C'(\xi)-C'(\xi)|=0.
\end{align*}

$\bullet$ Suppose the result holds true for some $n+1\leq N$.
By the induction hypothesis, $v'_{n+1}$ is in $BLP(\X')$, thus $R'v_{n+1}'$ is also in $BL(\X')$ by Proposition \ref{Rf-BL}.  For all $d\in\A$, $c'_d$ is in $BLP(\X')$ by Lemma \ref{lem-C'}.  As $\K(\xi=(\theta,\gamma))=\K(\gamma)$ depends only on $\gamma$ and $BLP(\X')$ is clearly stable by finite maximum, it follows that $v'_n$ is also in $BLP(\X')$.\\
Set  $\xi$ in $\bar\X'$. 
If $\xi=\Delta$, recall that $\K'(\Delta)=\{\cd\}$, so that on the one hand, one has
\begin{align*}
v'_n(\xi)&=c'(\xi,\cd)+R'v'_{n+1}(\xi,\cd)=0+v'_{n+1}(\Delta),
\end{align*}
and on the other hand, one has
\begin{align*}
\bar v'_n(\xi)&=\bar c'(\xi,\cd)+\bar R'\bar v'_{n+1}(\xi,\cd)=0+\bar v'_{n+1}(\Delta).
\end{align*}
Hence, one has $|v'_n(\Delta)-\bar v'_n(\Delta)|=|v'_{n+1}(\Delta)-\bar v'_{n+1}(\Delta)|=0$
by induction.\\
If $\xi=(\thetabar,\gamma)\in\bar \X'-\{\Delta\}$, with $\gamma=(y,w)$, then by definition one has $\theta(
\Omega_{m_0})\geq \big({\underline f}{\overline f}^{-1}\big)^{\frac{w}{\delta}\vee 1}e^{-w\|\lambda\|}$. The dynamic programming equations yield
\begin{align*}
|v'_n(\thetabar,\gamma)-\bar v'_n(\thetabar,\gamma)|
&\leq  \max_{d\in\K'(\xi)}|c'(\xi,d)-\bar c'(\xi,d)|+\max_{d\in\K'(\xi)}|R'v'_{n+1}(\xi,d)-\bar R'v'_{n+1}(\xi,d)|\\
&\quad +\max_{d\in\K'(\xi)}|\bar R'v'_{n+1}(\xi,d)-\bar R'\bar v'_{n+1}(\xi,d)|
\quad =A_1+A_2+A_3.
\end{align*}

The fist term $A_1$ is bounded thanks to Lemma \ref{lem-c'barc'} by $(B_c+L_c) \sup_{j \in \{1,\dots,\ell\}} \mathcal{D}_j$.

The second term $A_2$ is bounded by Proposition \ref{prop:R-Rbar} as $v'_n$ is in $BLP(\X')$
\begin{align*}
A_2 
&\leq (\|v'_{n+1}\|+[v'_{n+1}]){\bar C_{R'}} \big({\underline f}^{-1}{\overline f}\big)^{N\vee 1}e^{N\delta\|\lambda\|}\sup_{j \in \{1,\dots,\ell\}} \mathcal{D}_j.
\end{align*}

The last term $A_3$ is bounded by the induction hypothesis and using the fact that $\bar R'$ is a Markov kernel
\begin{align*}
A_3&\leq \sup_{\xi'\in \bar \X'}|v'_{n+1}(\xi')-\bar v'_{n+1}(\xi')|\leq \|{v'_{n+1}}\| \sup_{j \in \{1,\dots,\ell\}} \mathcal{D}_j,
\end{align*}
hence the result.
%
\subsection{Error bounds for the second discretization}
\label{app:error2}
%
The proof of Theorem \ref{theorem2} follows similar lines as that of Theorem \ref{theorem1}. Again, they are based on Lipschitz regularity properties of the operators involved, thus proofs are omitted. 
We introduce additional function spaces. 
Let $BLP(\bar\X')$ be the set of real-valued, bounded measurable functions $\varphi$ on $\bar\X'$ that are Lipschitz continuous and set
\begin{align*}
\|\varphi\|=\sup_{\xi\in \bar\X'}\|\varphi(\xi)\|,\quad
[\varphi]=\sup_{(\theta,\gamma)\neq (\theta',\gamma)\in \bar\X'-\{\Delta\}}\frac{|\varphi(\theta,\gamma)-\varphi(\thetabar ,\gamma)|}{d_E(\theta,\thetabar )}.
\end{align*}
Denote also $BLP_1(\bar\X')  =  \big\{ \varphi \in BLP(\bar\X') : \|\varphi\|_{E,\mathcal{P}}+ [\varphi]_{E,\mathcal{P}} \leq 1  \big\}$.
%
%
%
\subsubsection{Regularity of operator $\bar R'$} 
%
We first need to investigate the regularity of the approximate filter operator.
%
\begin{lem}\label{Psibar-theta-thetabar}
Under Assumption \ref{H1}, there exists some positive constant $C_\Psibar$ such that for all $\gamma'=(y',w')\in\O$, $d=(\ell,r) \in L\times\T$ and  all $\theta,\thetabar \in \mathcal{P}(\Omega)$ satisfying 
$\theta(\Omega_{m_0})>0$, $\thetabar(\Omega_{m_0})>0$,
 we have
\begin{align*}
\int_I d_\Omega(\Psibar(\theta,y',z',d),\Psibar(\thetabar,y',z',d))dy' \leq\frac{ C_\Psibar}{\theta(\Omega_{m_0})\thetabar(\Omega_{m_0})} d_\Omega(\theta,\thetabar).
\end{align*}
\end{lem}
%
Now we can turn to the regularity of operator $\bar R'$.
%
\begin{prop}\label{Rfbar-BL}
Under Assumption \ref{H1}, there exists a positive constant $C_{\bar R'}$ such that for all $g \in BLP_1(\bar\X')$, for all $\theta, \thetabar $ in $\mathcal P(\Omega)$, $(y,\bar y,w) \in I^2\times\{0,1\}\times[0,H]$ and $d=(\ell,r) \in \A-\{\cd\}$, such that $(\theta,y,w)\in\bar X'$ and $(\thetabar,\bar y,w)\in\bar X'$, 
one has
\begin{align*}
\big|\bar R'g(\theta,y,z,d)-\bar R'g(\thetabar ,\bar y,z,d)\big|
& \leq C_{\bar R'} d_\Omega(\theta,\thetabar ).
\end{align*}
In particular, operator $\bar R'$ maps $BLP(\bar\X')$ onto itself and for $g\in BLP(\bar\X')$, one has $\|\bar R'g\|\leq \|g\|$ and $[\bar R'g]\leq (\|g\|+[g])C_{\bar R'}$.
\end{prop}
%
\begin{prop}\label{prop:Rbar-Rhat}
For all $g \in BLP_1(\bar\X')$ and $1\leq j\leq n_{\Gamma}$, one has 
\begin{align*}
\left|\bar R'g(\rho^j,d)-\hat R'g(\rho^j,d)\right| \leq \sup_{\{1\leq j\leq n_{\Gamma}\}} \mathcal{\bar D}_j.
\end{align*}
\end{prop}
%
\begin{proof}
One has
\begin{align*}
\left|\bar R'g(\rho^j,d)-\hat R'g(\rho^j,d)\right|
&= \left|\int_{\bar \X'}g\big(\xi\big) \bar R'(d\xi| \rho^j,d)-\sum_{k=1}^{n_{\Gamma}}g\big(\rho^k\big)\bar R'(\bar C_k|\rho^j,d)\right|\\
&\leq \sum_{k=1}^{n_{\Gamma}}  \int_{\bar C_{k}} \big| g(\xi) -g(\rho^k) \big|\bar R'(d\xi| \rho^j,d)\\
&\leq \sum_{k=1}^{n_{\Gamma}}  \int_{\bar C_{k}}d(\xi,\rho^k)\bar R'(d\xi| \rho^j,d)\ \leq \sup_{\{1\leq k\leq {n_{\Gamma}}\}} \mathcal{\bar D}_k,
\end{align*}
hence the result.
\end{proof}
\subsubsection{Regularity of the cost functions}
\label{app:proof-cC'}
%
We last need to check that the cost functions $c'$, $C'$ also belong to $BLP(\bar\X')$.
%
\begin{lem}\label{lem-C'Om}
Let $g$ be a function from $\X'$ onto $\R$ belonging to $BLP(\X)$, then the restriction of $g$ to $\bar X'$ is in $BLP(\bar\X')$.
\end{lem}
%
In particular, the restrictions of $C'$ and $\bar c'$ to $\bar X'$ belong to $BLP(\bar\X')$.
%
\subsubsection{Proof of Theorem \ref{theorem2}}
\label{app:proof-th2}
%
We establish the result by (backward) induction on $k$, with the additional statement that $\bar v'_k$ is in $BLP(\bar\X')$ for all $k$.

$\bullet$ For $k=N$ and $1\leq j\leq n_{\Gamma}$, one has
\begin{align*}
\left|\vbar'_N(\rho^j)-\hat v'(\rho^j)\right| = \left|C'(\rho^j)-C'(\rho^j) \right| = 0,
\end{align*}
by definition. In addition, $\vbar'=C'$ is in $BLP(\bar\X')$.

$\bullet$ Suppose the result holds true for some $k+1\leq N$.
By the induction hypothesis, $\bar v'_{n+1}$ is in $BLP(\bar\X')$, thus $\bar R'v_{n+1}'$ is also in $BLP(\bar\X')$ by Proposition \ref{Rfbar-BL}.  For all $d\in\A$, $\bar c'_d$ is in $BLP(\bar\X')$.  As $BLP(\bar\X')$ is clearly stable by finite maximum, it follows that $\bar v'_n$ is also in $BLP(\bar\X')$.
Now, the dynamic programming equations yield
\begin{align*}
  \left|\vbar'_k(\rho^j)-\hat v'_k(\rho^j)\right| 
& \leq  \max_{d \in \K'(\rho^j)}  \left|\bar R'\vbar'_{k+1}(\rho^j,d)-\hat R'\bar v'_{k+1}(\rho^j,d) \right| \\
&\cutline  + \max_{d \in \K'(\rho^j)}   \left|\hat R'\bar v'_{k+1}(\rho^j,d) -\hat R'\hat v'_{k+1}(\rho^j,d) \right| .
\end{align*}
The last term on the right-hand side is smaller than $\max_{1\leq j\leq n_{\Gamma}}|\vbar'_{k+1}(\rho^j)-\hat v'_{k+1}(\rho^j)|$ as $\hat R'$ is a Markov kernel.
The first term on the right-hand side is bounded by Proposition \ref{prop:Rbar-Rhat} as $\bar v'_n$ is in $BLP(\bar\X')$
\begin{align*}
{ \max_{d \in \K'(\rho^j)} |\bar R'\bar v'_{k+1}(\rho^j,d)-\hat R'\bar v'_{k+1}(\rho^j,d) | }
&\leq C_{\bar R'} (\|\bar v'_{k+1}\| +[\bar v'_{k+1}])\sup_{j \in \{1,\dots,{n_{\Gamma}}\}} \mathcal{\bar D}_j.
\end{align*}
which concludes the proof.
%
\section{Additional details for the simulation study}
\label{supp:local}
%
The simulation study presented in \cref{sec:simulation} of the main manuscript has been constructed based on real data obtained from the \textit{Centre de Recherche en Cancérologie de Toulouse} (CRCT). Multiple myeloma (MM) is the second most common haematological malignancy in the world and is characterized by the accumulation of malignant plasma cells in the bone marrow. Classical treatments are based on chemotherapies, which, if appropriate, act fast and efficiently bringing MM patients to remission in a few weeks. However almost all patients eventually relapse more than once and the five-year survival rate is about $50\%$. 

We have obtained data from the \textit{Intergroupe Francophone du Myélome} 2009 clinical trial which has followed $748$ French MM patients from diagnosis to their first relapse on a standardized protocol for up to six years. At each visit a blood sample has been obtained to evaluate the amount of monoclonal immunoglobulin protein in the blood, a marker for the disease progression. 
Based on these data, we chose to use exponential flows for $\bs\Phi$, piece-wise constant linear functions for jump intensities, and three possible visit values: $\T=\{\delta, 2\delta, 4\delta\}$ with $\delta=15$ days. Explicit forms are given in \cref{supp:dynamics}, assumptions are verified in \cref{app:tech-spec}.

\subsection{Local characteristics in the simulation study}
\label{supp:local}
%

We now detail the special form used in our numerical examples to fit with our medical decision problem and prove that our assumptions are satisfied. We choose $\zeta_0=1$, $D=40$, $H=2400$ days and $\T=\{15,30,60\}$ days. 

\subsubsection{Special form of the local characteristics}
\label{supp:dynamics}

The values of $\Phi_m^\ell$ are given in \cref{tab:flot}, where $v_1^\emptyset=0.02$, $v_2^\emptyset=0.006$, $v'_1=0.077$, $v'_2=0.025$, $v_1=0.01$ and $v_2=0.003$. 
The link function $F$ is chosen to be the identity, and the noise $\varepsilon$ corresponds to a truncated centred Gaussian noise with variance parameter $\sigma^2$ and truncation parameter $s$. The explicit form of $t^{*d}_m$ is given in \cref{tab:tstar}.
\begin{table}[htp]
\footnotesize
\caption{Flow of the controlled PDMP}
\label{tab:flot}
\centering
\scalebox{0.8}{
\begin{tabular}{l|lll}
&$\ell=\emptyset$&$\ell=a$&$\ell=b$\\
\hline
$m=0$&$\Phi_m^\ell(\zeta_0,t)=\zeta_0$&$\Phi_m^\ell(\zeta_0,t)=\zeta_0$&$\Phi_m^\ell(\zeta_0,t)=\zeta_0$\\
$m=1$&$\Phi_m^\ell(\zeta,t)=\zeta e^{v_1^\emptyset t}$&$\Phi_m^\ell(\zeta,t)=\zeta e^{-v'_1 t}$&$\Phi_m^\ell(\zeta,t)=\zeta e^{v_1t}$\\
$m=2$&$\Phi_m^\ell(\zeta,t)=\zeta e^{v_2^\emptyset t}$&$\Phi_m^\ell(\zeta,t)=\zeta e^{v_2t}$&$\Phi_m^\ell(\zeta,t)=\zeta e^{-v'_2t}$\\
$m=3$&$\Phi_m^\ell(D,t)=D$&$\Phi_m^\ell(D,t)=D$&$\Phi_m^\ell(D,t)=D$
\end{tabular}}
\end{table}
%
%
\begin{table}[htp]
\footnotesize
\caption{Deterministic time to reach the boundary of the state space for the controlled PDMP.}
\label{tab:tstar} 
\centering
\scalebox{0.85}{
\begin{tabular}{l|lll}
$m$&$\ell=\emptyset$&$\ell=a$&$\ell=b$\\
\hline
$0$&$t^{*d}_m(\zeta_0)=+\infty$&$t^{*d}_m(\zeta_0)=+\infty$&$t^{*d}_m(\zeta_0)=+\infty$\\
$1$&$t^{*d}_m(\zeta)=\frac{1}{v_1^\emptyset}\log \frac{D}{\zeta}$&$t^{*d}_m(\zeta)=\frac{1}{v_1'}\log \frac{\zeta}{\zeta_0}$&$t^{*d}_m(\zeta)=\frac{1}{v_1}\log \frac{D}{\zeta}$\\
$2$&$t^{*d}_m(\zeta)=\frac{1}{v_2^\emptyset}\log \frac{D}{\zeta}$&$t^{*d}_m(\zeta)=\frac{1}{v_2}\log \frac{D}{\zeta}$&$t^{*d}_m(\zeta)=\frac{1}{v_2'}\log \frac{\zeta}{\zeta_0}$\\
\end{tabular}}
\end{table}

The values of $\lambda_m^\ell$ are given in \cref{tab:intensite}. For the standard relapse intensities ($\mu_i$), we choose  piece-wise increasing linear functions calibrated such that the risk of relapsing increases until $\tau_1$ (average of standard relapses occurrences), then remains constant, and further increases between $\tau_2$ and $\tau_3$ years (to model  late or non-relapsing patients):
\begin{align*}
\mu_i(x)=\left[ \begin{array}{ll}
\frac{\nu_1^i}{\tau_1^i}x & \text{ if } x\in[0,\tau_1^i] \\ 
\nu^i_1 & \text{ if } x\in[\tau_1^i,\tau_2] \\ 
\nu^i_1+\frac{\nu^i_2-\nu^i_1}{\tau_3-\tau_2}(x-\tau_2) & \text{ if } x\in[\tau_2,\tau_3] \\ 
\nu^i_2 & \text{ if } x \geq \tau_3
\end{array}\right.
\end{align*}
We set $\tau_1^1=750$, $\tau_1^2=500$ (days), $\tau_2=5$, $\tau_3=6$ (years), $\nu^i_1$ was selected so that $20$\% of patients relapse before $\tau^i_1$, and $\nu^i_2$ such that $10$\% of patients have not relapsed at horizon time $H$. 

For the therapeutic escape relapses (patients who relapse while treated for a current relapse), we chose to fit a Weibull survival distribution of the form
\begin{align*}
\mu_i'(\zeta)&=(\beta_i'\zeta)^{\alpha_i'},
\end{align*}
with $-1<\alpha_i'<0$ to account for a higher relapse risk when the marker decreases. We arbitrarily chose $\beta_i'=-0.8$ and calibrated $b_i=1000$ such that only about 5\% of patients experience a therapeutic escape.
\begin{table}[htp]
\footnotesize
\caption{Intensity of the controlled continuous time PDMP.}
\label{tab:intensite}
\centering
\scalebox{0.8}{
\begin{tabular}{l|lll}
&$\ell=\emptyset$&$\ell=a$&$\ell=b$\\
\hline
$m=0$&$\lambda_m^\ell(\zeta,u)=(\mu_1\!+\!\mu_2)(u)$&$\lambda_m^\ell(\zeta,u)=\mu_2(u)$&$\lambda_m^\ell(\zeta,u)=\mu_1(u)$\\
$m=1$&$\lambda_m^\ell(\zeta,u)=0$&$\lambda_m^\ell(\zeta,u)=\mu'_2(\zeta)$&$\lambda_m^\ell(\zeta,u)=0$\\
$m=2$&$\lambda_m^\ell(\zeta,u)=0$&$\lambda_m^\ell(\zeta,u)=0$&$\lambda_m^\ell(\zeta,u)=\mu'_1(\zeta)$ \\
$m=3$&$\lambda_m^\ell(D,u)=0$&$\lambda_m^\ell(D,u)=0$&$\lambda_m^\ell(D,u)=0$
\end{tabular}}
\end{table}

Finally, the Markov kernels are given in \cref{tab:Q}. Cases for $m=3$ are omitted as no jump is allowed when the patient has died.
The possible transitions between modes are illustrated in \cref{fig:states}, and an example of(continuous-time) controlled trajectory is given in \cref{fig:dynamics}.
\begin{table*}[htp]
\caption{\label{tab:Q} Markov kernel of the controlled continuous time PDMP, for $(A,B)$ Borel subsets of $(\zeta_0,D)$, and $[0,H]$ respectively}
\footnotesize
\centering
{
\begin{tabular}{l|l}
&$\ell=\emptyset$\\
\hline
$m=0$&$\bs Q^\ell(\{m'\}\times A\times B| m,\zeta_0,u)=\1_A(\zeta_0)\1_B(0)\1_{(m'\in\{1,2\})}\frac{\mu_{m'}(u)}{\mu_1(u)+\mu_2(u)}$\\
$m=1$&$\bs Q^\ell(\{m'\}\times A\times B| m,\zeta,u)=\1_A(\zeta)\1_B(0)\1_{(m'=0)}\1_{(\zeta=\zeta_0)}$\\
	   &$\bs Q^\ell(\{3\}\times A\times B| m,\zeta,u)=\1_A(D)\1_B(0))\1_{(\zeta=D)}$\\
$m=2$&$\bs Q^\ell(\{m'\}\times A\times B| m,\zeta,u)=\1_A(\zeta)\1_B(0))\1_{(m'=0)}\1_{(\zeta=\zeta_0)}$\\
	   &$\bs Q^\ell(\{3\}\times A\times B| m,\zeta,u)=\1_A(D)\1_B(0)\1_{(\zeta=D)}$\\
\hline
\hline
&$\ell=a$\\
\hline
$m=0$&$\bs Q^\ell(\{m'\}\times A\times B| m,\zeta_0,u)=\1_A(\zeta_0)\1_B(0)\1_{(m'=2)}$\\
$m=1$ &$\bs Q^\ell(\{m'\}\times A\times B| m,\zeta,u)=\1_A(\zeta)\1_B(0)\1_{(m'=2)}\1_{(\zeta>\zeta_0)}$\\
	   &$\bs Q^\ell(\{0\}\times A\times B| m,\zeta,u)=\1_A(\zeta_0)\1_B(0)\1_{(\zeta=\zeta_0)}$\\
$m=2$&$\bs Q^\ell(\{0,1\}\times A\times B| m,\zeta,u)=0$\\
	   &$\bs Q^\ell(\{3\}\times A\times B| m,\zeta,u)=\1_A(D)\1_B(0)\1_{(\zeta=D)}$\\
\hline
\hline
&$\ell=b$\\
\hline
$m=0$&$\bs Q^\ell(\{m'\}\times A\times B| m,\zeta_0,u)=\1_A(\zeta_0)\1_B(0)\1_{(m'=1)}$\\
$m=1$&$\bs Q^\ell(\{0,2\}\times A\times B| m,\zeta,u)=0$\\
	   &$\bs Q^\ell(\{3\}\times A\times B| m,\zeta,u)=\1_A(D)\1_B(0)\1_{(\zeta=D)}$\\
$m=2$&$\bs Q^\ell(\{m'\}\times A\times B| m,\zeta,u)=\1_A(\zeta)\1_B(0)\1_{(m'=1)}\1_{(\zeta>\zeta_0)}$\\
	   &$\bs Q^\ell(\{0\}\times A\times B| m,\zeta,u)=\1_A(\zeta_0)\1_B(0))\1_{(\zeta=\zeta_0)}$\\
\end{tabular}}
\end{table*}
%
\begin{figure}[htp]
\begin{center}
\begin{tikzpicture}
\node[draw,rounded corners=2pt] (mode0) at (0,0) {$0$};
\node[draw,rounded corners=2pt] (mode1) at (1.5,1) {$1$};
\node[draw,rounded corners=2pt] (mode2) at (1.5,-1) {$2$};
\node[draw,rounded corners=2pt] (cemetery) at (3,0) {$3$};

\draw[->,>=latex,dashed] (mode2.north east) -- (mode1.south east) node[right=0.01cm,pos=0.5]{\tiny{$b$}};
\draw[->,>=latex,dashed] (mode1.south west) -- (mode2.north west) node[left=0.01cm,pos=0.5]{\tiny{$a$}};
\draw[->,>=latex] (mode2) -- (mode0) node[right=0.01cm,pos=0.7]{\tiny{$b$}};
\draw[->,>=latex] (mode1) -- (cemetery) node[right=0.01cm,pos=0.3]{\tiny{$\neq a$}};
\draw[->,>=latex] (mode2) -- (cemetery) node[right=0.01cm,pos=0.3]{\tiny{$\neq b$}};
\draw[->,>=latex] (mode1) -- (mode0) node[right=0.01cm,pos=0.7]{\tiny{$a$}};
\draw[->,>=latex,dashed] (mode0.south west) -- (mode2.south west) node[left=0.01cm,pos=0.6]{\tiny{$\neq b$}};
\draw[->,>=latex,dashed] (mode0.north west) -- (mode1.north west) node[left=0.01cm,pos=0.6]{\tiny{$\neq a$}};
\end{tikzpicture}
\end{center}
\caption{\label{fig:states} {\bf State graph.} Full line arrows indicate deterministic jumps at the boundary, while dashed arrows indicate stochastic jumps. Letters indicate under which treatments the jumps are possible.}
\end{figure}

\begin{figure}[htp]
\begin{center}
      \begin{tikzpicture}[,scale=0.65, every node/.style={scale=0.65}]
      \coordinate (t0) at (-2, 0);
      \coordinate(T1) at (-0.2, 0) ;
      \coordinate (d1) at (1.5,1);
      \coordinate (T2) at (2.6,0);
      \coordinate (d2) at (2.8,0);
      \coordinate (T3) at (3.2,0);
      \coordinate (d3) at (4.5,1.5);
      \coordinate (d4) at (5,1.8);
      \coordinate (T4) at (6.1,0);
      \coordinate (d5) at (6.5,0);
      \coordinate (H) at (8,0);
      \draw[-] (t0)--(T1) 
      node[above=0.1,pos=0.1]{$\zeta_0,\small{\emptyset}$} node[below=0.1,pos=1] {$T_1$};
      \draw[color=green!60] (T1) to[out=0,in=-80] (d1);
      \draw[color=green!60] (d1) to[out=-90,in=-180] (T2) ;
      \node at (1.5,1.3) {\small{$a$}};
      \draw[-] (T2)--(d2) node[below=0.1,pos=-0.1,color=black] {$T_2$};
      \draw[-] (d2)--(T3) node[below=0.1,pos=1,color=black] {$T_3$}; 
      \node at (3,0.4) {\small{$\emptyset$}};
      \draw[color=red!60] (T3) to[out=0,in=-90] (d3);
      \draw[color=red!60] (d3) to[out=0,in=-100] (d4);
      \node at (4.5,1.75) {\small{$a$}};
      \draw[color=red!60] (d4) to[out=-90,in=-180] (T4);
      \node at (5,2.1) {\small{$b$}};
      \draw[-] (T4)--(d5) node[below=0.1,pos=-0.1,color=black] {$T_4$};
      \draw[-] (d5)--(H) node[below=0.1,pos=1,color=black] {$H$};
      \draw[-,color=black!50](-2,-0.1) -- (-2,0.1);
      \draw[-,color=black!50](-1,-0.1) -- (-1,0.1);
      \draw[-,color=black!50](0,-0.1) -- (0,0.1);
      \draw[-,color=black!50](1,0.15) -- (1,0.35);
      \draw[-,color=black!50](1.5,0.8) -- (1.5,1);
      \draw[-,color=black!50](2,0.1) -- (2,0.3);
      \draw[-,color=black!50](2.5,-0.1) -- (2.5,0.1);
      \draw[-,color=black!50](3,-0.1) -- (3,0.1);
      \draw[-,color=black!50](3.5,-0.05) -- (3.5,0.15);
      \draw[-,color=black!50](4.5,1.4) -- (4.5,1.6);
      \draw[-,color=black!50](5,1.7) -- (5,1.9);
      \draw[-,color=black!50](5.5,0.15) -- (5.5,0.35);
      \draw[-,color=black!50](6,-0.1) -- (6,0.1);
      \draw[-,color=black!50](7,-0.1) -- (7,0.1);
      \draw[-,color=black!50](7.5,-0.1) -- (7.5,0.1);
      \draw[-,color=black!50](8,-0.1) -- (8,0.1);
      \draw[dashed] (-2,3)--(8,3) node[above=0.15cm,pos=0] {$D$};
      \node at (7,0.4) {\small{$\emptyset$}};
    \end{tikzpicture}
\end{center}
\caption{\label{fig:dynamics} {\bf Example of controlled trajectory.} Representation of coordinate $\zeta$ of a controlled process. The process is in mode $0$ when drawn in black, in mode $1$ when in green, and in mode $2$ when in red. The $T_i$ represent hidden natural jump times. Letters above indicate the changes of treatment that can only occur at dates in the time grid $\delta^{1:N}=\{\delta,2\delta,\ldots, N\delta\}$.} 
\end{figure}

For the observations, we use the identity link function and a truncated gaussian for the noise. We arbitrarily chose $\sigma^2=1$ for the variance of the Gaussian distribution and  $s=2$ for the truncation parameter. 
%
\subsubsection{Technical specifications}
\label{app:tech-spec}

Assumption \ref{H-Phi} is valid for our exponential flows with $[\Phi]\leq e^{\bar{v}\delta}$ and $[t^{*}]\leq \frac{1}{\zeta_0\underline{v}}$ with $\underline{v}=\min\{v_1^\emptyset,v_1,v_2,v_2^\emptyset,v'_1,v'_2\}$.
In addition, as $\Phi_m^d(\zeta,t)\in E$ for $t< t^{*d}_m(\zeta)$ one also has  $|\Phi_m^d(\zeta,t)|\leq D$ and if $m\neq 0$, $|t^{*d}_m(\zeta)|\leq\frac{1}{\underline{v}}\log \frac{D}{\zeta_0}$. 

Assumptions \ref{H-lambda-lip}  and \ref{H-lambda-b} are valid for our intensities. 
For all $x\in E$, one has 
$$\|\lambda\|\leq \max\{\nu_2^1,\nu_2^2, (\beta_1'\zeta_0)^{\alpha_1'},(\beta_2'\zeta_0)^{\alpha_2'}\}.$$
The Lipschitz property is valid for $\mu'_i$ as $\zeta_0>0$ and $\alpha'_i>-1$. It is also valid for $\mu_i$ with $[\lambda]\leq \max\{\nu_1^i/\tau_1^i;\nu_2^i-\nu_1^i/\tau_3-\tau_2\}$.


Assumption \ref{H1} is valid for our truncated Gaussian noise and identity link function  with $L_f=2sD(D+s)(p\sigma^3\sqrt{2\pi})^{-1}$, $\bar f\leq (p\sigma\sqrt{2\pi})^{-1}$, $\underline{f}\geq (p\sigma\sqrt{2\pi})^{-1}e^{-D^2/2\sigma^2}$ and $B_f=2(D+s)\bar f$ where $p=\P(-s\leq Z\leq s)$ for a centred Gaussian random variable $Z$ with variance $\sigma^2$.

Assumption \ref{Hc} is valid for both the time-dependent and marker-dependent cost functions with  parameters given in \cref{tab:cost-reg}.
\begin{table*}[htp]
\footnotesize
\caption{\label{tab:cost-reg} Upper bounds for the regularity parameter for the cost functions.}
\centering
\scalebox{1}{
\begin{tabular}{r|l|l}
&time-dependent cost & marker-dependent cost \\ \hline
$B_C\leq$ & $\max\{|D-\zeta_0|,c_M\}$ & $\max\{|D-\zeta_0|,c_M\}$ \\
$B_c\leq$ &$ \max\{\max \T \times\max\{\gamma_i,\beta_i\},c_M\}$& $\max\{\max \T \times(\max\{\gamma_i,\beta_i\},|D-\zeta_0|),c_M\}$\\
$L_C \leq$ &$ 2 |D-\zeta_0|$& $ 2 |D-\zeta_0|$\\
$L_c \leq$& $2\max\T \max\{\gamma_i,\beta_i\}$& $2 \max \T (|D-\zeta_0|+\max\{\gamma_i,\beta_i\})$
\end{tabular} }
\end{table*}
%
\subsection{Further discussion on parameter tuning}
\label{supp:tuning}
%
\subsubsection{Grids construction}
\label{supp:grids}
%
As explained in the main manuscript, starting from an initial grid (with $184=n_\Omega$ points chosen to emphasize strong beliefs on each atom of grid $\Omega$), grids were extended iteratively by simulating a number $n_{sim}$ of trajectories using optimal strategies obtained by dynamic programming on the previous grids, and including all estimated filters with distance to their projection larger than a fixed threshold $s$. 

We varied the couples ($(n_{sim}, s)$ for both distances ($L_2$ and $L_m$) at each iteration. 
Extensive simulations (work not shown) indicates that for a similar number of points in the resulting grids, the results are better when using a small number of simulations with a stringent threshold than a larger number of simulations with a larger threshold. 

At each iteration, we also removed from the current grid all points whose density did not reach a given threshold. To do so, we simulated $10000$ trajectories using the current optimal strategy keeping track of all visited points in the grid. Then we removed all points with density smaller than $0.001/n_\Gamma$ (note that if all points were used equally, the density would be $1/n_\Gamma$). Due to the large amount of points removed at the first iteration, the second grid was computed without additional points.
The final process is illustrated in \cref{fig:supp-grids}.

\begin{figure}[htp]
\begin{center}
\tikzstyle{best} = [rectangle,font=\ttfamily, draw,text centered, text height=1.5ex, text depth=.25ex, minimum size=2.5ex,color=red]
\tikzstyle{intbest} = [rectangle,font=\ttfamily, draw,text centered, text height=1.5ex, text depth=.25ex, minimum size=2.5ex,color=green]
 \begin{tikzpicture}[scale=0.6]
   \node  at (0,2) {~}; 
   \node (init) at (2,0) {\textcolor{blue}{184}}; 
   \node  at (-3.5,0) {$\Gamma^0$}; 
   \node (g1) at (2,-2) {\textcolor{blue}{319}}; 
   \node  at (-3.5,-2) {$\Gamma^1$}; 
   \draw[->] (init.south)  -- (g1.north) node[left=0.1cm,pos=0.3] {\footnotesize $(10,0.2)$} node[left=0.35,pos=0.7]{\footnotesize $+93$} node[right=0.1cm,pos=0.7] {\footnotesize $-42$};
   \node (g2) at (2,-4.5) {\textcolor{blue}{549}}; 
   \node  at (-3.5,-4.5) {$\Gamma^2$}; 
   \draw[->] (g1.south)  -- (g2.north) node[left=0.1cm,pos=0.3] {\footnotesize $(10,0.2)$} node[left=0.35,pos=0.7]{\footnotesize $+205$} node[right=0.1cm,pos=0.7] {\footnotesize $-25$};
   \node (g3) at (2,-7) {\textcolor{blue}{722}}; 
   \node  at (-3.5,-7) {$\Gamma^3$}; 
   \draw[->] (g2.south)  -- (g3.north) node[left=0.1cm,pos=0.3] {\footnotesize $(100,0.4)$} node[left=0.35,pos=0.7]{\footnotesize $+163$} node[right=0.1cm,pos=0.7] {\footnotesize $-10$};
   \node (g4) at (2,-9.5) {\textcolor{blue}{989}}; 
   \draw[->] (g3.south)  -- (g4.north) node[left=0.1cm,pos=0.3] {\footnotesize $(12,0.2)$} node[left=0.35,pos=0.7]{\footnotesize $+266$} node[right=0.1cm,pos=0.7] {\footnotesize $-1$};
   \node  at (-3.5,-9.5) {$\Gamma^4$}; 
   \node (initm) at (10,0) {\textcolor{blue}{184}}; 
   \node (g1m) at (10,-2) {\textcolor{blue}{310}}; 
   \draw[->] (initm.south)  -- (g1m.north) node[left=0.1cm,pos=0.3] {\footnotesize $(7,0.2)$} node[left=0.35,pos=0.7]{\footnotesize $+81$} node[right=0.1cm,pos=0.7] {\footnotesize $-45$};
   \node (g2m) at (10,-4.5) {\textcolor{blue}{545}}; 
   \draw[->] (g1m.south)  -- (g2m.north) node[left=0.1cm,pos=0.3] {\footnotesize $(7,0.2)$} node[left=0.35,pos=0.7]{\footnotesize $+213$} node[right=0.1cm,pos=0.7] {\footnotesize $-22$};
   \node (g3m) at (10,-7) {\textcolor{blue}{709}}; 
   \draw[->] (g2m.south)  -- (g3m.north) node[left=0.1cm,pos=0.3] {\footnotesize $(40,0.55)$} node[left=0.35,pos=0.7]{\footnotesize $+160$} node[right=0.1cm,pos=0.7] {\footnotesize $-4$};
   \node (g4m) at (10,-9.5) {\textcolor{blue}{1021}}; 
   \draw[->] (g3m.south)  -- (g4m.north) node[left=0.1cm,pos=0.3] {\footnotesize $(8,0.2)$} node[left=0.35,pos=0.7]{\footnotesize $+307$} node[right=0.1cm,pos=0.7] {\footnotesize $-5$};
 \end{tikzpicture}
 \caption{\label{fig:supp-grids} Iterative grids for distance $L_2$ (left) and $L_m$ (right)}
\end{center}
\end{figure}
%
\subsubsection{Distance impact on trajectory}
\label{app:traj}
%
Setting the same seed, we illustrate the impact of the grid and distance choice on a trajectory in \cref{fig:L2,fig:Lm}. The first row shows the (true) value of $X_2$, with $X_1$ indicated by circles for $X_1=0$, triangles for $X_2=1$ and pluses for $X_2=3$. The second row shows the observed process $Y$, with colors indicating treatments: black for $\emptyset$, green for $a$ and red for $b$. The third row shows the mass probability of each mode of the estimated filter, and the fourth for its projected counterpart. Finally, the fifth row shows the distance between estimated and projected filters, in black for the $L_2$ distance, and in blue for the $L_m$ distance. 

Figure \ref{fig:L2} shows that though the distance between the filter and its projection significantly decreases between the two grids, the projection does not preserve the mode, hence leads to decision which do not always seem appropriate seeing the data, in particular regarding the next visit date. On the contrary, \cref{fig:Lm} shows that distance $L_m$ decreases the distance meanwhile maintaining the mass distribution between modes.
\begin{figure}[htp]
\centering
\includegraphics[width=\textwidth]{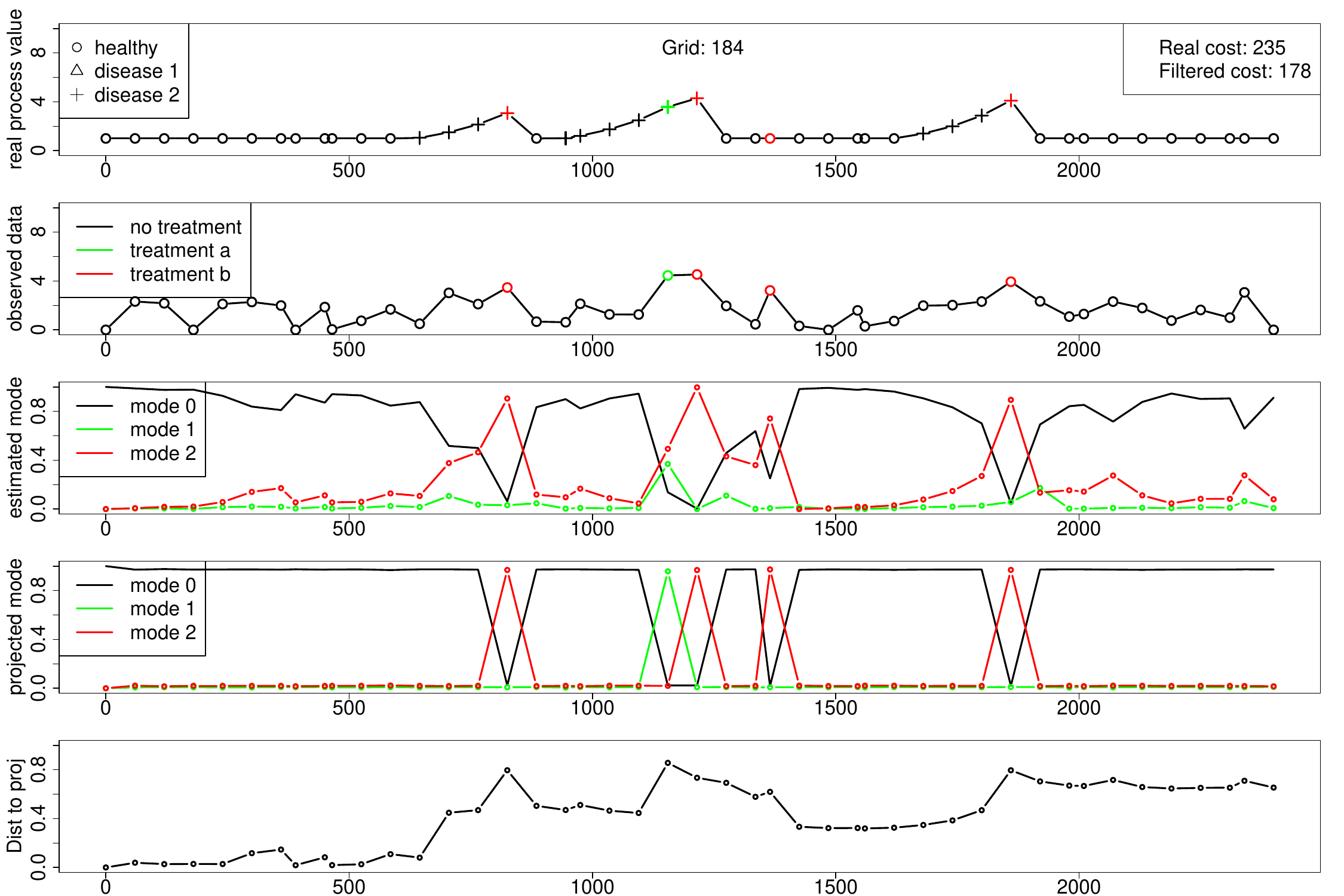} \\
\vspace{1cm}
\includegraphics[width=\textwidth]{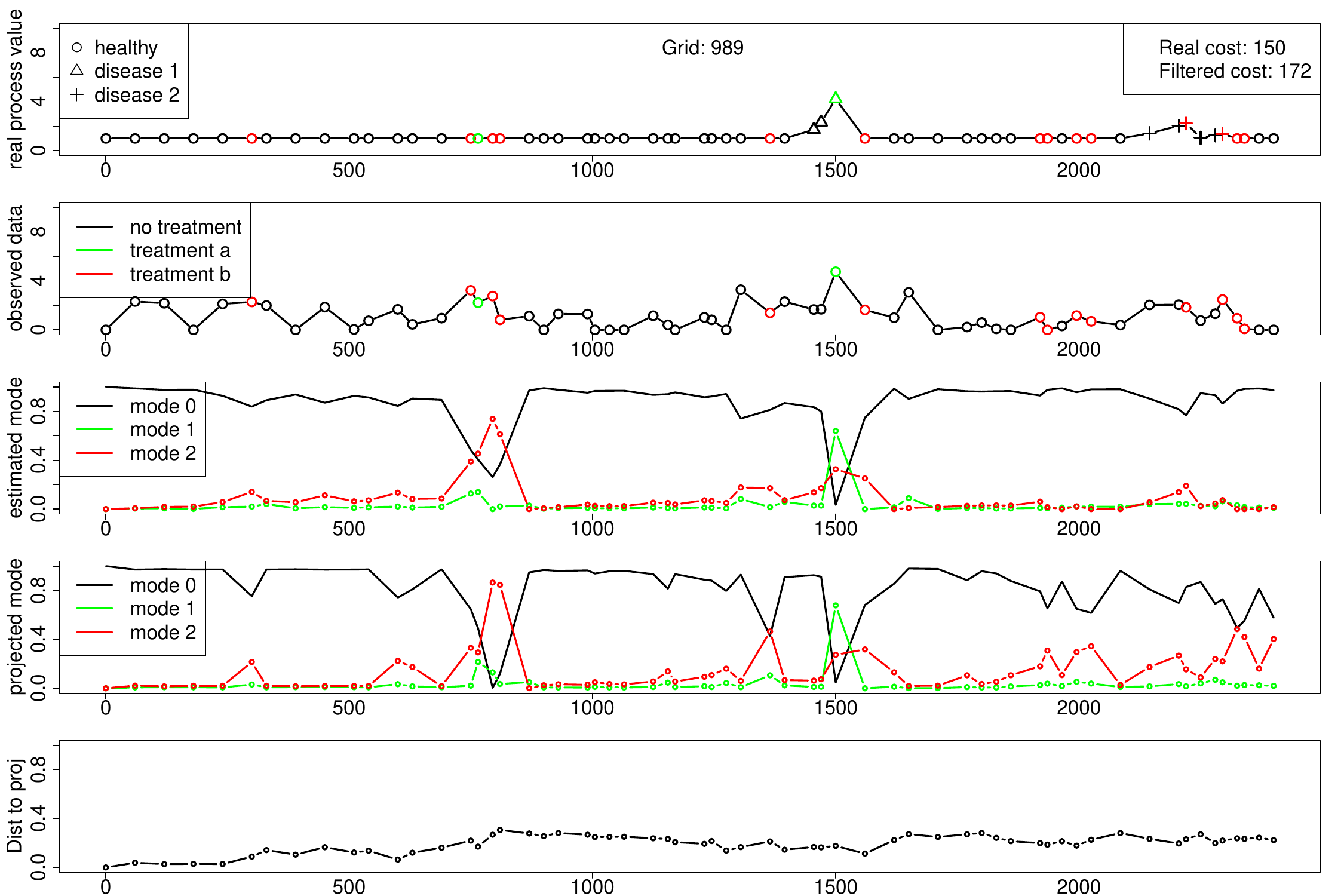} 
\caption{\label{fig:L2} {\bf Trajectories for $L_2$ distance.}}
\end{figure}
\begin{figure}[htp]
\centering
\includegraphics[width=\textwidth]{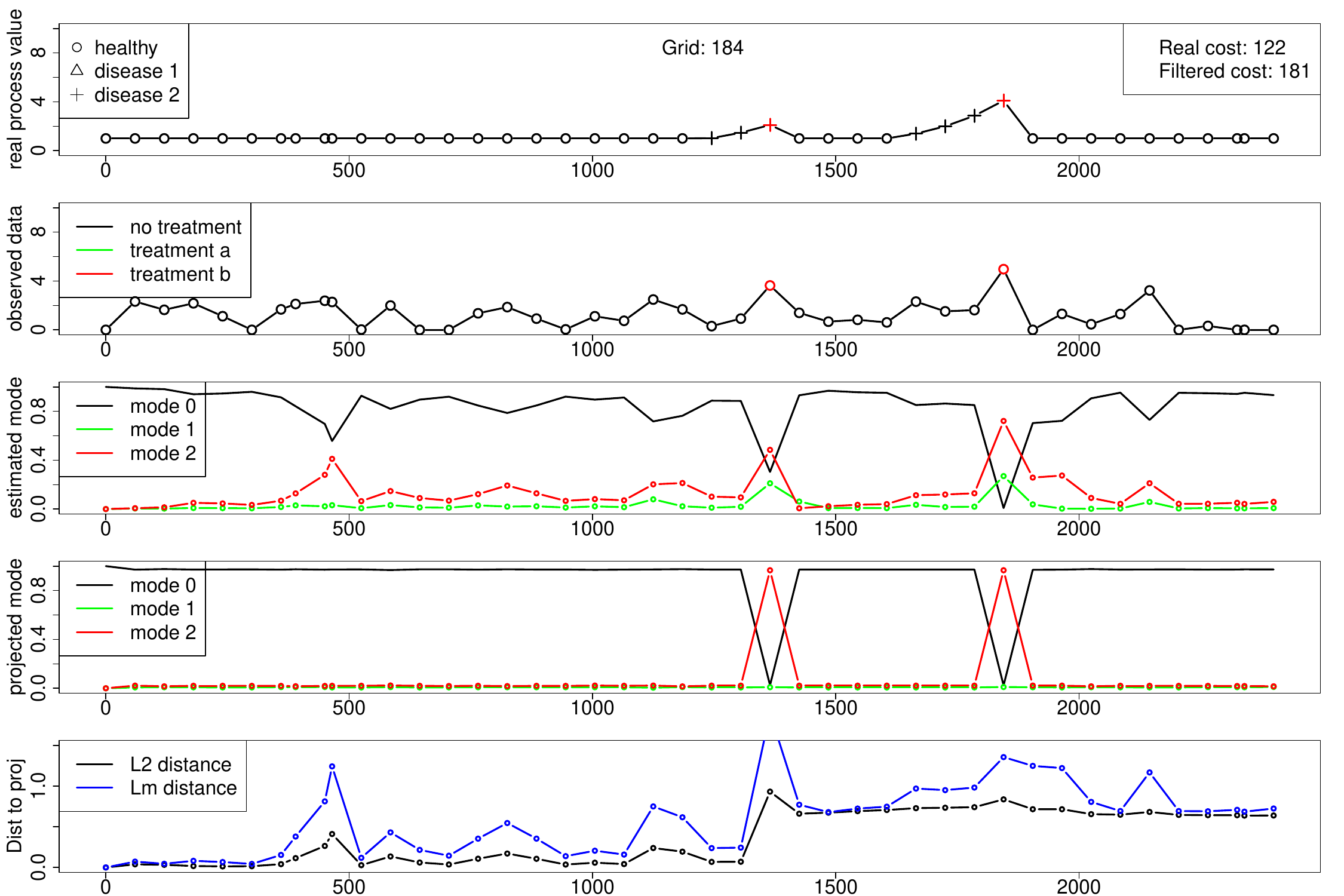} \\
\vspace{1cm}
\includegraphics[width=\textwidth]{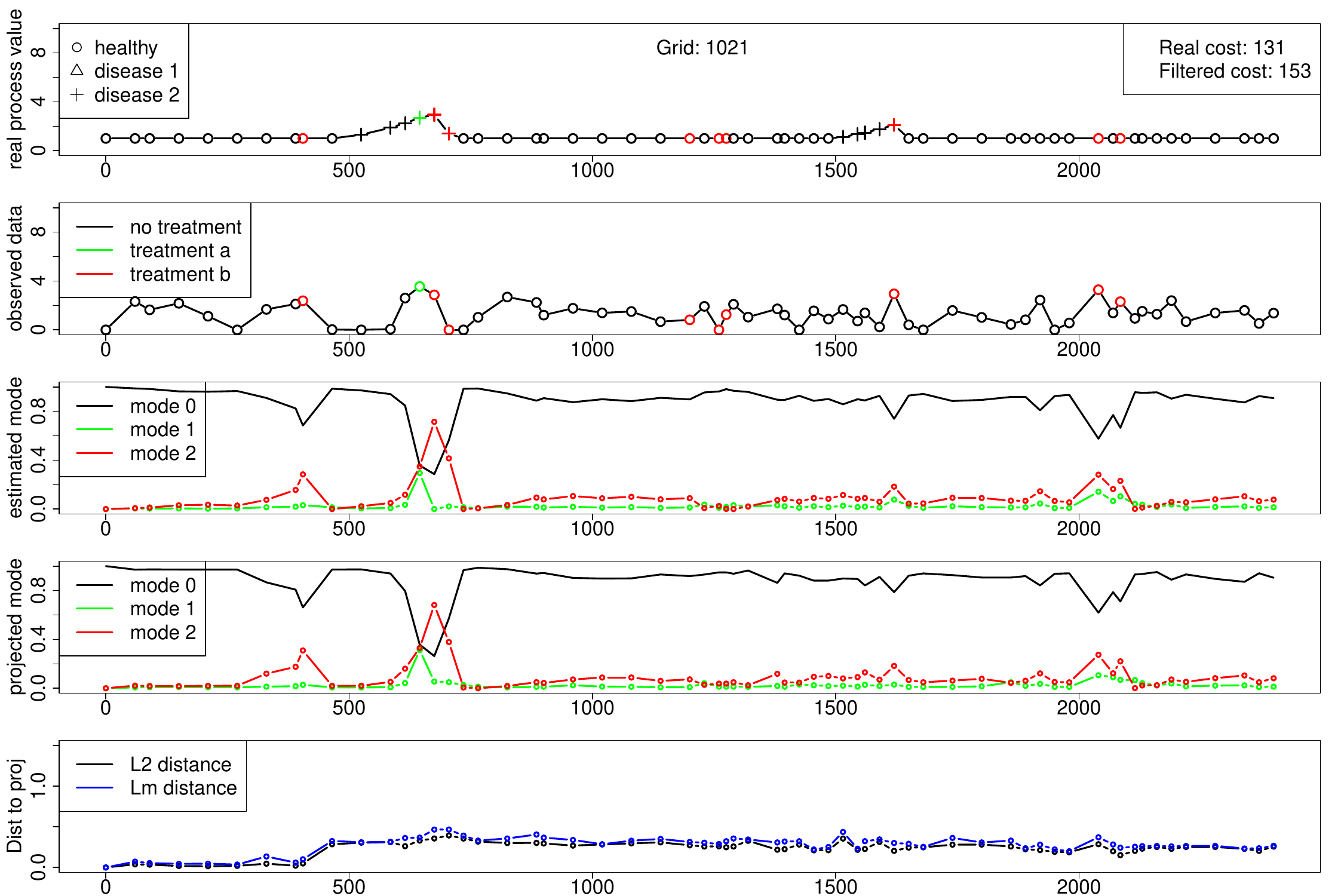} 
\caption{\label{fig:Lm} {\bf Trajectories for $L_m$ distance.}}
\end{figure}
%
\subsubsection{Choice of $\hat \Psi$ or $\bar \Psi$ in practice}
\label{app:operator}
%
The discretization strategy lead to a Markov kernel $\hat{R}'$ on $\bar\X'$ to which can be associated a Markov chain $\hat\Psi$. In theory, the dynamic programming algorithm operates on this Markov chain, and error bounds from the main theorem are computed accordingly. This should imply that at iteration $k$, filter $\bar\Psi_k$ is computed from the new observation and the current filter $\hat\Psi_{k-1}$, then projected on $\Gamma$ to identify the optimal decision. Then this projection $\hat\Psi_k=p_{\Gamma}(\bar\Psi_k)$ is saved as current filter for the next observation. 

In practice, this implies that the projection error $\hat\Psi_k-\bar\Psi_k$ is propagated through the dynamic programming recursion. We propose instead to save $\bar\Psi_k$ as current filter for the next iteration, \textit{i.e.} the filter at iteration $k+1$ will be computed using $Y_{k+1}$ and $\bar\Psi_k$. We do not propose any error bound on this practical strategy, and there is no guarantee that this should lead to better results, as projection errors on $\Gamma$ may sometimes compensate previous errors. However in our simulation studies we have observed a significant difference between the strategies, as illustrated in \cref{tab:psi} (grid $1021$ with $L_m$ distance). \\
One can note that the estimated cost is identical in the visit choice framework between filter  $\bar \Psi$ and $\hat \Psi$ as they correspond to  how grids were calibrated, but  in practice the real cost is lower with $\bar \Psi$ as the latter is closer to the true data.

\begin{table}[htp]
\footnotesize
\caption{\label{tab:psi} Choice of $\hat \Psi$ or $\bar \Psi$ as current filter}
\centering
\begin{tabular}{c|r|ccc||cc}
& & & $\bar \Psi$ &  & $\hat \Psi$ & \\
\hline
 & Visits & $\hat{v}_0$ & real cost (sd) & est. cost (sd) &  real cost  (sd) & est. cost  (sd) \\ 
  \hline
 \hline
&OS& 132.43 & 136.96  (3.97) & 134.8  (0.75) & 161.38  (7.2) & 134.82  (0.96) \\ 
Grid&FD-15 & 253.28 & 214.29  (1.67) & 215.83  (0.78) & 227.84  (3.13) & 264.19  (2.31) \\ 
1231&FD-60 & 163.57 & 143.02  (5.29) & 141.8  (1.05) & 163.19  (7.29) & 165.22  (1.05) \\ 
   \hline
\end{tabular}
\end{table}


\begin{thebibliography}{10}

\bibitem{Alm01}
{\sc Almudevar, A.}
\newblock A dynamic programming algorithm for the optimal control of piecewise
  deterministic markov processes.
\newblock {\em SIAM Journal on Control and Optimization 40}, 2 (2001),
  525--539.

\bibitem{BL17}
{\sc B\"auerle, N., and Lange, D.}
\newblock Optimal control of partially observable piecewise deterministic
  {M}arkov processes.
\newblock {\em SIAM J. Control Optim. 56}, 2 (2018), 1441--1462.

\bibitem{bauerle11}
{\sc B\"auerle, N., and Rieder, U.}
\newblock {\em Markov decision processes with applications to finance}.
\newblock Universitext. Springer, Heidelberg, 2011.

\bibitem{baysse2013maintenance}
{\sc Baysse, C., Bihannic, D., G{\'e}gout-Petit, A., Prenat, M., de~Saporta,
  B., and Saracco, J.}
\newblock Maintenance optimisation of optronic equipment.
\newblock {\em Chemical Engineering Transactions 33\/} (2013), 709--714.

\bibitem{BBGPPS14}
{\sc Baysse, C., Bihannic, D., G\'egout-Petit, A., Prenat, M., and Saracco, J.}
\newblock Hidden {M}arkov model for the detection of a degraded operating mode
  of optronic equipment.
\newblock {\em J. SFdS 155}, 3 (2014), 48--61.

\bibitem{BdSD13}
{\sc Brandejsky, A., de~Saporta, B., and Dufour, F.}
\newblock Optimal stopping for partially observed piecewise-deterministic
  {M}arkov processes.
\newblock {\em Stochastic Process. Appl. 123}, 8 (2013), 3201--3238.

\bibitem{C20}
{\sc Calvia, A.}
\newblock Stochastic filtering and optimal control of pure jump markov
  processes with noise-free partial observation.
\newblock {\em ESAIM: Control, Optimisation and Calculus of Variations 26\/}
  (2020), 25.

\bibitem{Aut18}
{\sc Cleynen, A., and {d}e {S}aporta, B.}
\newblock Change-point detection for piecewise deterministic markov processes.
\newblock {\em Automatica 97\/} (2018), 234--247.

\bibitem{cocozza2021markov}
{\sc Cocozza-Thivent, C., et~al.}
\newblock {\em Markov renewal and piecewise deterministic processes}.
\newblock Springer, 2021.

\bibitem{costa_piecewise_2016}
{\sc Costa, M.}
\newblock A piecewise deterministic model for a prey-predator community.
\newblock {\em The Annals of Applied Probability 26}, 6 (Dec. 2016).

\bibitem{CD13}
{\sc Costa, O., and Dufour, F.}
\newblock {\em Continuous average control of piecewise deterministic {M}arkov
  processes}.
\newblock SpringerBriefs in Mathematics. Springer, New York, 2013.

\bibitem{CR00}
{\sc Costa, O., and Raymundo, C.}
\newblock Impulse and continuous control of piecewise deterministic markov
  processes.
\newblock {\em Stochastics and Stochastic Reports 70}, 1-2 (2000), 75--107.

\bibitem{CD89}
{\sc Costa, O.~L., and Davis, M.}
\newblock Impulse control of piecewise-deterministic processes.
\newblock {\em Mathematics of Control, Signals and Systems 2\/} (1989),
  187--206.

\bibitem{Davis84}
{\sc Davis, M.}
\newblock Piecewise-deterministic {M}arkov processes: a general class of
  nondiffusion stochastic models.
\newblock {\em J. Roy. Statist. Soc. Ser. B 46}, 3 (1984), 353--388.
\newblock With discussion.

\bibitem{Davis93}
{\sc Davis, M.}
\newblock {\em Markov models and optimization}, vol.~49 of {\em Monographs on
  Statistics and Applied Probability}.
\newblock Chapman \& Hall, London, 1993.

\bibitem{Aut12}
{\sc de~{S}aporta, B., and {D}ufour, F.}
\newblock Numerical method for impulse control of piecewise deterministic
  {M}arkov processes.
\newblock {\em Automatica 48}, 5 (2012), 779--793.

\bibitem{Aut16}
{\sc de~{S}aporta, B., Dufour, F., and Geeraert, A.}
\newblock Optimal strategies for impulse control of piecewise deterministic
  markov processes.
\newblock {\em Automatica 77\/} (2017), 219--229.

\bibitem{dSDZ16}
{\sc de~Saporta, B., Dufour, F., and Zhang, H.}
\newblock {\em Numerical methods for simulation and optimization of piecewise
  deterministic {M}arkov processes}.
\newblock Mathematics and Statistics Series. ISTE, London; John Wiley \& Sons,
  Inc., Hoboken, NJ, 2016.
\newblock Application to reliability.

\bibitem{de2019dynamic}
{\sc de~Saporta, B., Dufour, F., and Zhang, H.}
\newblock Dynamic optimization of maintenance policies for multistate system.
\newblock In {\em ESREL 2019-29th European Safety and Reliability Conference\/}
  (2019).

\bibitem{Dempster95}
{\sc Dempster, M., and Ye, J.}
\newblock Impulse control of piecewise deterministic {M}arkov processes.
\newblock {\em Ann. Appl. Probab. 5}, 2 (1995), 399--423.

\bibitem{devooght1997dynamic}
{\sc Devooght, J.}
\newblock Dynamic reliability.
\newblock In {\em Advances in nuclear science and technology}. Springer, 1997,
  pp.~215--278.

\bibitem{devooght97}
{\sc Devooght, J.}
\newblock {\em Dynamic reliability}.
\newblock Advances in nuclear science and technology. Chapman and Hall, Berlin,
  1997.

\bibitem{DHKL15}
{\sc Doumic, M., Hoffmann, M.and~Krell, N., and Robert, L.}
\newblock Statistical estimation of a growth-fragmentation model observed on a
  genealogical tree.
\newblock {\em Bernoulli 21}, 3 (2015), 1760--1799.

\bibitem{DHP16}
{\sc Dufour, F., Horiguchi, M., and Piunovskiy, A.}
\newblock Optimal impulsive control of piecewise deterministic markov
  processes.
\newblock {\em Stochastics 88}, 7 (2016), 1073--1098.

\bibitem{fritsch2015modeling}
{\sc Fritsch, C., Harmand, J., and Campillo, F.}
\newblock A modeling approach of the chemostat.
\newblock {\em Ecological Modelling 299\/} (2015), 1--13.

\bibitem{guo2009continuous}
{\sc Guo, X., Hern{\'a}ndez-Lerma, O., Guo, X., and Hern{\'a}ndez-Lerma, O.}
\newblock {\em Continuous-time Markov decision processes}.
\newblock Springer, 2009.

\bibitem{herbach_inferring_2017}
{\sc Herbach, U., Bonnaffoux, A., Espinasse, T., and Gandrillon, O.}
\newblock Inferring gene regulatory networks from single-cell data: a
  mechanistic approach.
\newblock {\em BMC Systems Biology 11}, 1 (Dec. 2017).

\bibitem{lange2017cost}
{\sc Lange, D.~K.}
\newblock {\em Cost optimal control of Piecewise Deterministic Markov Processes
  under partial observation}.
\newblock PhD thesis, Dissertation, Karlsruhe, Karlsruher Institut f{\"u}r
  Technologie (KIT), 2017, 2017.

\bibitem{makis2003optimal}
{\sc Makis, V., and Jiang, X.}
\newblock Optimal replacement under partial observations.
\newblock {\em Mathematics of operations research 28}, 2 (2003), 382--394.

\bibitem{MZ96}
{\sc Marseguerra, M., and Zio, E.}
\newblock {M}onte {C}arlo approach to {PSA} for dynamic process systems.
\newblock {\em Reliability Engineering and System Safety 52}, 3 SPEC. ISS.
  (1996), 227--241.

\bibitem{pasin_controlling_2018}
{\sc Pasin, C., Dufour, F., Villain, L., Zhang, H., and Thiébaut, R.}
\newblock Controlling {IL}-7 {Injections} in {HIV}-{Infected} {Patients}.
\newblock {\em Bulletin of Mathematical Biology 80}, 9 (Sept. 2018),
  2349--2377.

\bibitem{PPTZ19}
{\sc Piunovskiy, A., Plakhov, A., Torres, D.~F., and Zhang, Y.}
\newblock Optimal impulse control of dynamical systems.
\newblock {\em SIAM Journal on Control and Optimization 57}, 4 (2019),
  2720--2752.

\bibitem{RT15}
{\sc Riedler, M., and Thieullen, M.}
\newblock Spatio-temporal hybrid ({PDMP}) models: central limit theorem and
  {L}angevin approximation for global fluctuations. {A}pplication to
  electrophysiology.
\newblock {\em Bernoulli 21}, 2 (2015), 647--696.

\bibitem{RTW12}
{\sc Riedler, M., Thieullen, M., and Wainrib, G.}
\newblock Limit theorems for infinite-dimensional piecewise deterministic
  {M}arkov processes. {A}pplications to stochastic excitable membrane models.
\newblock {\em Electron. J. Probab. 17\/} (2012), no. 55, 48.

\bibitem{soleimani2021integration}
{\sc Soleimani, M., Campean, F., and Neagu, D.}
\newblock Integration of hidden markov modelling and bayesian network for fault
  detection and prediction of complex engineered systems.
\newblock {\em Reliability Engineering \& System Safety 215\/} (2021), 107808.

\bibitem{ZDDG08}
{\sc {Z}hang, H., {D}ufour, F., Dutuit, Y., and Gonzalez, K.}
\newblock Piecewise deterministic markov processes and dynamic reliability.
\newblock {\em Proc. of the Inst. of Mechanical Engineers, Part O: Journal of
  Risk and Reliability 222}, 4 (2008), 545--551.

\end{thebibliography}
\end{document}